\documentclass{amsart}
\usepackage{amsfonts, amsmath,amscd, amssymb, latexsym, mathrsfs, stmaryrd, verbatim, wasysym }
\usepackage{graphicx}
\usepackage[colorinlistoftodos]{todonotes}
\usepackage[english]{babel}
\usepackage[utf8x]{inputenc}
\usepackage[all]{xy}
\usepackage{amscd}

\usepackage{appendix}
\usepackage{amsthm}

\usepackage{appendix}
\usepackage{mathtools}

\usepackage{color}
\definecolor{darkgreen}{cmyk}{1,0,1,.2}
\definecolor{m}{rgb}{1,0.1,1}
\definecolor{green}{cmyk}{1,0,1,0}

\definecolor{test}{rgb}{1,0,0}
\definecolor{cmyk}{cmyk}{0,1,1,0}

\long\def\red#1{\textcolor {red}{#1}}

\long\def\green#1{\textcolor {green}{#1}}
\long\def\m#1{\textcolor {m}{#1}}

\newtheorem{Equation}{}[section]

\newtheorem{corollary}[Equation]{Corollary}
\newtheorem{definition}[Equation]{Definition}

\newtheorem{lemma}[Equation]{Lemma}
\newtheorem{proposition}[Equation]{Proposition}

\newtheorem{problem}[Equation]{Problem}
\newtheorem{remark}[Equation]{Remark}

\newtheorem{theorem}[Equation]{Theorem}

\def\Dir{\operatorname{{\not{\hspace{-0.05cm}\pa}}}}

\def\ch{\operatorname{ch}}
\def\Ch{\operatorname{Ch}}
\def\codim{\operatorname{codim}}
\def\Diff{\operatorname{Diff}}
\def\Mor{\operatorname{Mor}}
\def\Dom{\operatorname{Dom}}
\def\Aut{\operatorname{Aut}}
\def\SO{\operatorname{SO}}
\def\SU{\operatorname{SU}}

\def\Dir{\operatorname{{\not{\hspace{-0.05cm}\pa}}}}

\def\reg{\operatorname{reg}}
\def\av{\operatorname{av}}

\def\Dom{\operatorname{Dom}}

\def\ind{\operatorname{ind}}
\def\End{\operatorname{End}}

\def\geo{\operatorname{geo}}

\def\Id{\operatorname{Id}}

\def\Hom{\operatorname{Hom}}

\def\Im{\operatorname{Im}}
\def\Ind{\operatorname{Ind}}

\def\Ker{\operatorname{Ker}}

\def\mod{\operatorname{mod}}

\def\Supp{\operatorname{Supp}}

\def\Spin{\operatorname{Spin}}

\def\Sp{\operatorname{Sp}}

\def\Tr{\operatorname{Tr}}
\def\STr{\operatorname{STr}}
\def\tr{\operatorname{tr}}

\def\C{\mathbb C}
\def\D{\mathbb D}
\def\H{\mathbb H}
\def\N{\mathbb N}
\def\R{\mathbb R}
\def\S{\mathbb S}

\def\Z{\mathbb Z}
\def\B{\mathbb B}
\def\LL{\mathbb L}

\def\Q{\mathbb Q}

\def\cA{{\mathcal A}}
\def\maA{{\mathcal A}}
\def\maa{{\what{\mathfrak a}}}

\def\maB{{\mathcal B}}
\def\maC{{\mathcal C}}
\def\maD{{\mathcal D}}

\def\maE{{\mathcal E}}

\def\maF{{\mathcal F}}
\def\maM{{\mathcal M}}
\def\cG{{\mathcal G}}

\def\maG{{\mathcal G}}
\def\cG{{\mathcal G}}

\def\cH{{\mathcal H}}

\def\maK{{\mathcal K}}
\def\maN{{\mathcal N}}

\def\maR{{\mathcal R}}

\def\maS{{\mathcal S}}
\def\maQ{{\mathcal Q}}

\def\maU{{\mathcal U}}
\def\cU{{\mathcal U}}

\def\what{\widehat}
\def\wtit{\widetilde}

\def\tD{\widetilde{D}} 
\def\tDp{\widetilde{D'}} 
\def\tF{\widetilde{F}} 
\def\tPhi{\widetilde{\Phi}}

\def\tE{\widetilde{E}}
\def\tH{\widetilde{H}}
\def\tX{\widetilde{X}}

\def\tW{\widetilde{W}}
\def\tw{\widetilde{w}}

\def\tG{\widetilde{G}}
 
\def\tF{\widetilde{F}} 
\def\tO{\widetilde{\Omega}}

\def\tM{\widetilde{M}}

\def\tS{\widetilde{S}}

\def\tT{\widetilde{T}}

\def\tX{\widetilde{X}}
\def\tZ{\widetilde{Z}}
\def\tg{\widetilde{g}}
\def\tQ{\widetilde{Q}}

\def\tY{\widetilde{Y}}
\def\tW{\widetilde{W}}
\def\tN{\widetilde{N}}

\def\tm{\widetilde{m}}

\def\dd{\displaystyle}
\def\pa{\partial}

\def\ep{\epsilon}

\begin{document}



\title[Dirac operators and PSC] 
{The Gromov-Lawson index   and the \\Baum-Connes assembly map}


\author[M-T. Benameur]{Moulay-Tahar Benameur}
\address{Institut Montpellierain Alexander Grothendieck, UMR 5149 du CNRS, Universit\'e de Montpellier}
\email{moulay.benameur@umontpellier.fr}

\thanks{MSC (2010) 53C12, 57R30, 53C27, 32Q10. }

\begin{abstract}
This contribution gathers some results related with the study of Positive Scalar Curvature metrics  in connection with the Baum-Connes assembly map. 
 \end{abstract}
\maketitle
\tableofcontents

%


\section*{Introduction}

Spin-Dirac operators provide the main obstructions to the existence of metrics with Positive Scalar Curvature (PSC metrics), and when such metrics exist, they are also involved in the study of moduli spaces associated with PSC metrics. The first known obstruction is the $A$-hat genus, as implied by   the Atiyah-Singer index formula for the Dirac operator $\Dir$, and by the Lichnerowicz computation of the spin Laplacian in terms of the scalar curvature in \cite{Lichne}. The argument is standard: The Lichnerowicz computation relates the spin-Dirac operator with the scalar curvature function $\kappa$:
$$
\Dir^2 \,\, = \,\, \nabla^* \nabla  +   \frac{\kappa}{4} \Id.
$$
Hence the index of $\Dir$ vanishes under a PSC metric. The Atiyah-Singer formula \cite{AtiyahSinger1, AtiyahSinger3} shows that the index of $\Dir$ is precisely the A-hat genus which hence should vanish under a PSC metric. 

Later in the work of Hitchin on harmonic spinors, the refined $\maa$-genus was used for dimensions which are not multiple of $4$ \cite{Hitchin}, and Lusztig showed how the higher genera involving the fundamental group (in the case of the free abelian group) also provide a larger set of obstructions \cite{Lusztig}. Using the spin-Dirac operator  now twisted by appropriate bundles, and involving more systematically the fundamental group of the manifold, Gromov and Lawson proved  some strong obstruction results. For instance,  they proved that enlargeable spin manifolds cannot hold PSC metrics, a byproduct of this theorem is that tori cannot hold PSC metrics \cite{GromovLawson2}, a longstanding question which  was  as well answered by Schoen and Yau  using a totally different approach based on their famous splitting theorem \cite{SchoenYau1, SchoenYau2}. The Gromov-Lawson approach through the Dirac operator led  them to the proof of an important  relative index theorem for pairs of complete non-compact  manifolds, and to also propose some higher invariants, see the seminal papers \cite{GromovLawson2, GromovLawson1, GromovLawson3}. Later on, Connes proved the first obstruction result for foliations in \cite{ConnesFundamental}. More precisely, he showed that for any closed (not necessarily spin) oriented manifold $M$ with $\what{A}(M)\neq 0$, no spin foliation $F$ on $M$ can hold PSC along the leaves. 
More recently, in \cite{Zhang}, Zhang obtained the same result but assuming $M$ to be spin and not necessarily the foliation. 

As observed by Rosenberg \cite{Rosenberg},  there is a close connection between these higher genera on a given spin manifold $M$ and the so-called Baum-Connes conjecture  for the fundamental group $\Gamma=\pi_1(M)$ of $M$. Recall that this conjecture states that an assembly map (BC map):
$$
\mu_\Gamma \, : \, K_*^{top} (\Gamma) \longrightarrow K_*(C_r^*\Gamma),
$$
 built up from the higher index construction, is an isomorphism.  In fact, the rational injectivity of $\mu_\Gamma$ implies the vanishing of the higher A-hat genera when $M$ has a PSC metric. Moreover, and at least for complete manifolds with cylindrical ends as built from compact manifold with boundary by the APS construction \cite{APS1}, The Gromov-Lawson index together with its $L^2$ version as introduced in \cite{BenameurCovering} are related with reduced eta invariants. Now the bijectivity of the (maximal) BC map  allows to prove many vanishing results for these reduced eta invariants under PSC metrics and when the fundamental group is torsion-free.  Indeed, when PSC metrics exist, while the index vanishes, still the spin-Dirac operator  provides these more subtle secondary  invariants, also called  rho invariants.  They  involve by definition the fundamental group through its finite dimensional unitary representations \cite{APS3} or more generally its regular $\ell^2$ representation through the so-called Cheeger-Gromov rho invariant \cite{CheegerGromov} which is defined as follows:
$$
\rho_{(2)} (\Dir) := \eta_{(2)} (\wtit\Dir) - \eta (\Dir) \text{ with }\wtit\Dir \text{ being the spin-Dirac operator on the universal covering}.
$$
This latter $L^2$ invariant   was efficiently used to distinguish non-concordant PSC metrics when $\Gamma$ has torsion. The relation of this invariant, and hence of the Gromov-Lawson index in this case, with the BC assembly map for the fundamental group $\Gamma$ was  discovered in the late 90's and has been progressively clarified, during the last two decades. When the fundamental group is torsion-free, the Cheeger-Gromov rho invariant of the spin-Dirac operator seems to behave like the spin-genus and to vanish under a PSC metric. This  was  proved only under the assumption that  the (maximal) BC conjecture is true, showing a strong implication of the BC assembly map. See for instance \cite{Keswani, PiazzaSchick2005, BenameurRoyJFA}. A kind of converse influence could also be progressively  observed, namely that  there is a potential impact of the index theory for PSC near infinity spin manifolds on the properties of the BC map, through the (higher) absolute index classes that the Gromov-Lawson construction provides, see \cite{XieYu2014}.

We have chosen in this survey paper to review the relation between the Gromov-Lawson (absolute and) relative index for complete spin manifolds with PSC near infinity and the BC conjecture. The rough slogan is that these constructions provide  tools for non-compact manifolds  similar, although different in nature for the absolute index,  to the ones provided by the Atiyah-Singer theorem for closed manifolds, at least when $\Gamma$ is torsion-free. In particular, one can, in principle, follow the known extensions  of the classical index theorem with their consequences, to  construct new invariants for moduli spaces of  PSC metrics. 
For this reason, we have added in Subsection \ref{BC-overview} below, for the convenience of the reader,  a brief overview of some generalizations of the classical Atiyah-Singer theorem, in connection with the BC assembly map.  \\

As for the Gromov-Lawson index theory on complete non-compact spin manifolds with PSC near infinity, some generalizations have  been carried out recently but the picture is still incomplete. For coverings,  the higher Gromov-Lawson  index was defined by Xie and Yu in  \cite{XieYu2014} where they also prove the higher relative index theorem, extending in this way the Fomenko-Miscenko approach. This higher index raised many interesting questions, and revealed a new potential connection with the BC map. The Atiyah $L^2$ index theory was also recently investigated in \cite{BenameurCovering} and the $L^2$-relative index theorem was deduced by using the Murray-von Neumann dimension theory. In Section \ref{SectionCovering}, we give a detailed survey of these developpements for coverings and also gather some related questions.  

Then comes the questions of generalizations to families, equivariants families and  foliations. The Connes measured index theorem has been recently extended in \cite{BenameurHeitsch21} into a  relative measured index theorem \`a la Gromov-Lawson for complete foliations with leafwise PSC near infinity. 

In Sections 2-4, some more or less new results are explained. 
The first one is the proof that leafwise PSC off some compact subspace of a spin complete foliation allows to define the higher spin index in the $K$-theory of the space of leaves, exactly as was proved by Connes for closed foliated manifold, completing the picture for foliations \cite{ConnesSurvey}. This is the exact generalization of a Gromov-Lawson result for complete spin manifolds. Notice that the trace of the compact subspace on a fixed leaf needs not to be compact in that leaf. 

\begin{theorem}
Let $(M, F)$ be a complete foliation with even dimensional leaves and let $D=D^\pm$ be a leafwise generalized Dirac operator which is invertible near infinity, then the operator $D^+$ has a well defined higher index class $\Ind (D^+)$ in the $K$-theory group $K(C^*(M, F))$ of the space of leaves. 
\end{theorem}

In fact we focus on the following important case:
\begin{corollary}
For the spin-Dirac operator twisted by a hermtian bundle with  connection $(E, \nabla^E)$, the spin-index class $\Spin (M, F; (g, \nabla^E))$ is well defined in $K(C^*(M, F))$, if we assume that  outside some compact subspace of $M$:
\begin{enumerate}
\item the leafwise scalar curvature of $g$ is uniformly positive;
\item  $\nabla^E$ is leafwise flat.
\end{enumerate}
\end{corollary}

An interesting situation is the case of $\Gamma$ actions on a compact manifold $X$, with $\Gamma=\pi_1(W)$  for a smooth even complete spin manifold $W$, and where one is interested in equivariant  families, parametrized by $X$, of metrics on the universal covering $\tW$, such that the family of scalar curvature functions, a global $\Gamma$-invariant smooth function on $\tW\times X$ is a uniformly positive function  outside some  $\Gamma$-compact subspace of $\tW\times X$. The previous covering situation corresponds to $X=\{\bullet\}$ but now we have a smooth non trivial foliation which is associated to this action by using the suspension construction. The higher index  lives here in the $K$-theory group of the crossed product $C^*$-algebra $C(X)\rtimes \Gamma$. The index and rho invariants in the classical closed case were investigated in \cite{BenameurPiazza}.

For a foliation given by a fibration $M\to B$ with even dimensional connected fibers and connected base for simplicity, we get  an index class in $K(B)$.  The interesting invariant $\iota (g_0, g_1)$ of Gromov-Lawson \cite{GromovLawson3} can then be generalized to families to provide an invariant in $H^{[0]} (B, \R)$. More precisely, given a fibration $N\to B$ now with odd spin fibers and with two metrics $g_0$ and $g_1$ having positive fiberwise scalar curvature, one defines as in \cite{GromovLawson3} the invariant $
\iota (N\vert B; g_0, g_1) \; \in \; H^{[0]} (B, \R)$,
to be the Chern character of the index of the fiberwise spin-Dirac operator on the fibration $M=N\times \R\to B$. When $E$ is an extra hermitian bundle over $N$ with fiberwise flat connections $\nabla^0$ and $\nabla^1$, then one gets similarly the invariant 
$$
\iota \left(N\vert B; (g_0, \nabla^0);  (g_1, \nabla^1)\right) \in  H^{[0]} (B, \Q).
$$
 So, whenever  $(g_0, \nabla^0)$ can be connected to $(g_1, \nabla^1)$ through positive fiberwise scalar curvature metrics and fiberwise flat connections, we see that $\iota \left(N\vert B; (g_0, \nabla^0);  (g_1, \nabla^1)\right)$ vanishes. This Gromov-Lawson  invariant satisfies the similar properties proved for a single manifold in \cite{GromovLawson3}. 

We also use the Bismut-Cheeger index formula for compact fibrations with boundary \cite{BismutCheeger1, MelrosePiazza1}, to outline a method to construct  a higher Kreck-Stolz  invariant $\sigma (N\vert B, g)$ for some smooth $4k-1$-fibrations with positive fiberwise scalar curvature under the right assumptions, especially  when the tangent bundle to the fibers has trivial real Pontrjagin classes. This invariant now lives in $H^{[0]} (B, \R)$ and as for the classical one, is an absolute version of the Gromov-Lawson invariant $\iota (N\vert B; g_0, g_1)$. This class is defined in terms of eta forms and the following is then a consequence of the Bismut-Cheeger formula.

\begin{proposition}
$\sigma (N, g)$ only depends on the fiberwise concordance class of $g$. It is also unchanged  under actions of fiberwise diffeomorphisms. 
\end{proposition}

In  Section \ref{RelativeFamily}, given two complete non-compact fibrations $M_0\to B_0$ and $M_1\to B_1$ over compact manifolds, with the two fiberwise generalized  Dirac operators $D_0$ and $D_1$ which agree near infinity, we prove, under usual assumptions,  that the higher relative index class $\Ind (D_0, D_1)$ can be defined as an element of the $K$-theory of the base $B_0$ once the usual identifications are used. Denoting by $\varphi_0:B_0\to B_1$ the induced diffeomorphism between the bases, we can state the formula for  the spin-Dirac operators with coefficients in bundles with hermitian connections which are also compatible near infinity:

\begin{theorem}\
The following formula holds in $H^{[0]} (B_0, \R)$:
$$
\Ch \Ind (\Dir_0\otimes \nabla^0, \Dir_1\otimes \nabla^1) = \int_{M_0\vert B_0} \what{A}(TM_0\vert B_0) \Ch (\nabla^0)\; - \;\varphi_0^*\int_{M_1\vert B_1} \what{A}(TM_1\vert B)\Ch (\nabla^1).
$$
\end{theorem}
The RHS is defined in a way similar to the case of a single manifold and makes sense. 


\ \\\subsection{From the AS-theorem to the BC-conjecture}\label{BC-overview}

We give a very brief overview of {\underline{some}} developpements that occured since the Atiyah-Singer index theorem was proved in the early 60's. See  \cite{AtiyahSinger63} and also \cite{AtiyahSinger1, AtiyahSinger3}.
There are nowadays many well  known geometric situations where the classical index theory has been fully explored and  efficiently used.  Already in the late 60's Atiyah and Singer extended their index formula to skew adjoint operators and real $K$-theory  \cite{AtiyahSinger5}  and also to families of closed manifolds \cite{AtiyahSinger4} which were as well stated in real $K$-theory \cite{Hitchin}. The  families index theorem was then efficiently used to prove the Novikov conjecture for the abelian free group \cite{Lusztig} and one deduced the vanishing of higher genera for the spin-Dirac operator under PSC by the same approach.  The  idea there was that the higher signatures, or higher genera, can be extracted from a  robust Atiyah-Singer  index bundle in the $K$-theory of the base  torus. By robust we mean here that the higher index class enjoys the allowed properties almost by construction because it benefits from the confortable properties of $C^*$-algebras. This means that it vanishes for instance for the family spin-Dirac operator under PSC along the fibers, and it is a fiberwise homotopy invariant for the family signature operator.  So, when the fundamental group is abelian, it has been  involved through  invariants of families parametrized by its Pontrjagin dual. When the group $\Gamma$ is non abelian, it needed to be replaced by its convolution $C^*$-algebra  $C^*_\Gamma$  which replaces the algebra of continuous functions on the Pontrjagin dual. 
Index theory for general (infinite) Galois  coverings started in the seminal work of Atiyah (and Singer)  \cite{AtiyahCovering}, giving the first implication of the Murray-von Neumann dimension theory, with a (deep) numerical  index equality: 
$$
\text{``index upstairs=index downstairs'',  }
$$
meaning that the $L^2$ index of the lifted operator to the covering coincides with the index of the operator on $M$. Fomenko and Miscenko later showed  in \cite{FomenkoMiscenko} that the Lusztig idea can be extended to all Galois coverings so as to provide an index formula in the $K$-theory of the $C^*$-algebra of the group, where the indices are again the robust classes. Unfortunately, despite the abelian case, the computation of the $K$-theory of the group $C^*$-algebra is a hard problem. The Fomenko-Miscenko construction  led though to the formulation of the famous Baum-Connes assembly map, expecting all the $K$-classes of the group $C^*$-algebra to be constructed in this way from an (generalized) index map. The Baum-Connes map $\mu_\Gamma$, alluded to above, is built up from the higher index construction,  organizing all these Dirac operators as a topological (in principle computable) group $K^{top}_*(\Gamma)$ by moding out with respect to equivalence relations preserving the higher index. In this way, the BC conjecture  proposes a topological method to compute the $K$-theory of the  group $C^*$-algebra, the  host for the robust index classes, so it confronts this hard problem. Another solution to get around the effective computation the $K$-theory of the group $C^*$-algebra is to try to extract numerical invariants using traces and also higher traces as developped by Connes through his cyclic cohomology \cite{ConnesIHES}, and one hopes to recover in this way the higher genera in the Gromov-Lawson-Rosenberg conjecture \cite{RosenbergGL} or the higher signatures in the Novikov conjecture \cite{Novikov}. This method succeded for Gromov-hyperbolic groups for instance \cite{ConnesMoscoviciHyper}, and also for a larger class of groups in \cite{ConnesGromovMoscovici}.

The next interesting extension of the Atiyah-Singer index theorem combines all the previous situations into the setting of foliations. Indeed, combining families and groups one is naturally led to group actions which on the other hand correspond through a standard suspension construction to a large class of interesting foliations. The case of group actions was also motivated by the need of a ``Baum-Connes conjecture with coefficients'' in the many attempts to prove the BC conjecture, and turned out to have interesting consequences in some solid physics problems \cite{Bellissard}. The first index theorem for foliations was proved by Connes in \cite{ConnesIntegration}, a far reaching generalization of Atiyah's $L^2$-index theorem which uses again the Murray-von Neumann dimension theory, but now in the context of groupoids. The higher index theorem for foliations was later proved by Connes and Skandalis \cite{ConnesSkandalis} following ideas from \cite{ConnesSurvey}; a higher equality in the $K$-theory of some $C^*$-algebra that Connes associated with the holonomy (or homotopy) groupoid of the foliation and which is the precise  replacement of the base in the case of fibrations, but also of the group $C^*$-algebra in the case of coverings when one uses the homotopy groupoid. A higher index theorem  for foliations was also proved later on, using the superconnection approach, in  \cite{HeitschLazarov} but now it is an equality in the Haefliger cohomology of the foliation. \\

{\em{Acknowledgements.}}  The author would like to thank J. Heitsch, V. Mathai and I. Roy  for several helpful discussions.


\section{BC assembly map and $L^2$ GL-index, an overview}\label{SectionCovering}

Let us start by briefly stating a natural  question which will be expanded in the sequel in connection with the BC assembly map and PSC metrics. Given a smooth complete manifold $M$ and a generalized Dirac operator $D$ which is invertible near infinity, in the sense that there exists $\kappa_0>0$ such that $D^2\geq \kappa_0$ off some compact subspace of $M$, it was proved in \cite{BenameurCovering} that the Atiyah $\Gamma$-index $\Ind_{(2)} (\tD)$ of the lift $\tD$ of $D$ to any Galois $\Gamma$-covering $\tM\to M$ is a well defined real number.

\begin{problem}\label{Integrality}
Can one find a pair $(M, D)$ as above with $\Gamma$ torsion-free and such that $\Ind_{(2)} (\tD)$ is not an integer? Can $M$ be chosen to be a complete spin manifold with PSC near infinity, such that $D$ is the spin-Dirac operator?
\end{problem}

As we shall see many other questions arise in close relation with Problem \ref{Integrality}. 

Recall that  a simply connected spin closed manifold  $M$ has a PSC metric  if and only if $\maa (M)=0$  \cite{Stolz}. For a non simply connected closed manifold $M$, the situation is much more complicated and one needs to involve appropriately the fundamental  group. An important observation is that the Lichnerowicz formula continues to hold for the spin-Dirac operator with coefficients in any flat bundle and even  when such flat bundle is a bundle of modules over some fixed unital $C^*$-algebra. There is a privileged flat bundle associated with the fundamental group called the Miscenko bundle whose  fibers are just the reduced $C^*$-algebra  $C_r^*\pi_1(M)$ viewed as a module over itself. The Lichnerowicz principle then gives an  obstruction to a PSC metric in terms of  an index of the spin-Dirac operator with coefficients in this Miscenko bundle,  a class in the $K$-theory group $K_{\dim (M)} (C_r^*\pi_1M)$ of $C^*_r\pi_1M$  that we shall define below in the more general case of PSC near infinity.

Let  now $\pi: \tM \to M$ be a Galois $\Gamma$-covering over the smooth complete  riemannian manifold $(M, g)$ which is not assumed to be compact. Here $\Gamma$ is any finitely generated countable discrete group. We endow $\tM$ with the lifted $\Gamma$-invariant metric $\tg$, so that it is also a complete riemannian manifold and $\pi$ is an isometric covering. We fix a hermitian bundle $S$ of generalized spinors over $M$, in the sense of \cite{GromovLawson3}[Section I]. Its pull-back to $\tM$ is denoted $\tS$. 
Assume that $M$ is even dimensional and that the generalized spin bundle $S$ has a $\Z_2$-grading $S=S^+\oplus S^-$, we then have the corresponding $\Gamma$-equivariant splitting   $\tS=\tS^+\oplus \tS^-$.  We fix corresponding  Lebesgue class measures on $M$ and $\tM$ that we denote by $dm$ or $d\tm$ respectively. 

Let $D$ be  a generalized Dirac operator  acting on the smooth sections of  $S$ over $M$, and denote by $\tD$ its $\Gamma$-invariant lift to $\tM$, acting on the sections of $\tS$:
$$
\tD: C_c^\infty (\tM, \tS) \longrightarrow C_c^\infty (\tM, \tS),\quad \tD = \left(\begin{array}{cc} 0 & \tD^-\\ \tD^+ & 0\end{array}\right).
$$
The $\Gamma$-action can be linearized to endow the space $C_c^\infty (\tM, \tS)$ with the structure of  a prehilbertian right module over the group convolution algebra $\C\Gamma$ of finitely supported complex functions on $\Gamma$.  More precisely, for $\xi, \xi'\in C_c^\infty (\tM, \tE)$ and $f\in \C\Gamma$, the rules are
$$
(\xi f) ( \tm) := \sum_{g\in\Gamma} (g\xi)(\tm) f(g^{-1}) \text{ and } < \xi , \xi' > (g) := \int_{\tM} <\xi (\tm), (g^{-1} \xi') (\tm)> d\tm,
$$
The completion of $\C\Gamma$ with respect to the regular representation in $\ell^2\Gamma$ is the regular $C^*$-algebra $C^*_r\Gamma$, and if we complete $C_c^\infty (\tM, \tE)$ with respect to the above inner product and the regular representation, we obtain a Hilbert module $\maE_r$ over $C^*_r\Gamma$. The Dirac operator $\tD$ then extends to a regular operator $\maD_r$ on $\maE_r$. \\

When the manifold $M$ is \underline{closed} and for instance even dimensional, the operator $D$ represents a $K$-homology class $[D]$  in $K_0 (M)$ which can be pushed forward as a class  $f_*[D]\in K_0(B\Gamma)$ under the classifying map $f:M\to B\Gamma$ of  the covering. On the other hand, the operator $\maD_r$ is also known to be  Fredholm and self-adjoint on $\maE_r$ so that its chiral part $\maD_r^+:\maE_r^+\to \maE_r^-$ has a well defined (higher) index class $\Ind (\maD_r)$ in  the $K$-theory group of the algebra $\maK (\maE_r)$ of compact operators on $\maE_r$. This is a standard construction of the so-called boundary map in $K$-theory. Now, the $K$-theory of the $C^*$-algebra $\maK (\maE_r)$ is just the $K$-theory of $C^*_r\Gamma$. Hence, we end up with a higher index class 
$$
\Ind (\maD_r) \in K_0(C^*_r\Gamma).
$$
It is a consequence again of the Lichnerowicz principle that if $M$ has PSC, then the class $\Ind (\maD_r)$ vanishes. The process of assigning to the spin (or more generally spin$^c$) {\underline{closed}} manifold $M$, the index class of the Dirac operator as above, including coefficients in hermitian bundles, is the building block for the  famous Baum-Connes assembly map (BC-map):
$$
\mu_\Gamma : K_*^{top} (\Gamma)\longrightarrow K_*(C_r^*\Gamma),
$$
where $K_*^{top} (\Gamma)$ is some (topological) representable $K$-homology group which reduces to the $K$-homology of the classifying space $B\Gamma$ when $\Gamma$ is torsion free. In general,  it takes into account  all proper (not necessarily free) cocompact $\Gamma$-actions. The BC conjecture states that the above map $\mu_\Gamma$ should be an isomorphism.  As a consequence of its very definition, any element of $K_*(C_r^*\Gamma)$ that is constructed as a higher index class of a Dirac operator as mentionned above  on closed manifolds, belongs to the range of the BC-map.  Completing with respect to all unitary representations of $\Gamma$, one gets the full $C^*$-algebra $C^*\Gamma$, and similarly a maximal Hilbert module $\maE$ over $C^*\Gamma$ and a maximal BC assembly map with the same domain but with the range replaced by the $K$-theory of the full algebra. For torsion-free groups, it is given as:
$$
K_*(B\Gamma) \longrightarrow K_*(C^*\Gamma).
$$
We say that $\Gamma$ satisfies the maximal BC conjecture if this latter assembly map is an isomorphism. Here is a short list of   examples where the BC conjecture is known.  If $\Gamma$ Gromov a-T-menable, especially if it is amenable, then the maximal (and reduced)  BC map is known to be an isomorphism. Any Gromov hyperbolic group satisfies the reduced BC conjecture, and the same is true for any cocompact discrete subgroup of $\Sp (n, 1)$ for $n>1$ and $SL(3, \R)$ by Lafforgue's thesis. Property (T) groups are not expected to satisfy the maximal BC conjecture since the natural group morphism $K_*(C_m^*\Gamma)\to K_*(C_r^*\Gamma)$ is not an isomorphism in this case.
\smallskip

Assume now that $M$ is a complete \underline{non-compact} even dimensional spin-manifold and that the complete metric $g$ has PSC outside some compact subspace of $M$. Then Gromov and Lawson proved that the Dirac operator $D$ has a well defined index $\Ind (D^+)\in \Z$, and introduced some concordance and diffeomorphism invariants of  PSC metrics on odd dimensional closed spin manifolds. Xie and Yu proved that for Galois $\Gamma$-coverings as above, the induced operator $\maD_r$ on the Hilbert module $\maE_r$ is again a regular  operator and induces a Fredholm operator with  a well defined higher index class $\Ind (\maD_r^+) \in K_0(C_r^*\Gamma)$ \cite{XieYu2014}. This allowed them  to  state and prove the higher relative index theorem of Gromov-Lawson. 
Now, for all groups for which the BC map is known to be surjective, the Xie-Yu higher index class must belong  to the range of the BC map.  Notice that the  $K$-theory index $\Ind (\maD_m)$ is in general not given in terms of the kernel and cokernel of the operator, and similarly for the reduced index $\Ind (\maD_r)$. These kernel modules do not define in general $K$-theory classes  without first perturbating appropriately the operator, see \cite{FomenkoMiscenko}.  

As observed by Rosenberg \cite{Rosenberg},  injectivity of the BC assembly map implies the vanishing of the higher $A$-hat genera of any spin manifold with PSC, a far-reaching generalization of the Lichnerowicz result, which is a consequence of the above vanishing of the higher index class. 

We come now to the relation with the $L^2$-index appearing in Problem \ref{Integrality}. In the closed case, applying the group morphism $K_0(C^*\Gamma)\to \Z$ induced by the trivial representation of $\Gamma$ in $\C$ (the average trace) to $\Ind (\maD_m)$,  one gets an integer which is nothing but the Fredholm index of $D^+$ on $M$. In the same way, it is well known in the closed case that applying to it the  group morphism $K_0(C^*\Gamma)\to \R$ induced by the regular trace (evaluation at the neutral element) one gets the $L^2$-index of $\tD$ on $\tM$ as introduced and studied by Atiyah in \cite{AtiyahCovering}.  The Atiyah $L^2$ index theorem \cite{AtiyahCovering}, say equality of these two scalar indices,  allows in particular to deduce that the $L^2$ index is an integer.  Rephrased  in terms of the BC assembly map in the torsion free case, the Atiyah theorem is the coincidence of the morphisms induced by these two traces  on the range of the maximal BC-map. Hence, (rational) surjectivity of $\mu_\Gamma$  implies in this case  the integrality of the $L^2$-index. But this latter integrality property is clearly true on any closed manifold without reference to the BC map. 

\begin{remark} When $\Gamma$ has  torsion,  these traces and others constructed out of the finite subgroups  provide  a measurement of some defect group \cite{WeinbergerYu}.   
\end{remark}

Back to complete non-compact manifolds and the relation between the Xie-Yu higher index class and the $L^2$-index, notice that the $L^2$-version of the Gromov-Lawson index that we have used in Problem \ref{Integrality} was defined  in terms of  kernel and cokernel thanks to the following almost-gap-in-the spectrum property:

\begin{theorem}\cite{BenameurCovering}\label{AlmostGap}\
Assume  that a generalized Dirac operator $D$ acting on the sections of $S$ over the complete riemannian manifold $M$ is uniformly invertible  near infinity. Then, there exists $\ep>0$ such that the spectral projection of $\tD$ for $(-\ep, +\ep)$ has finite $L^2$-dimension. 
\end{theorem}

In particular, the Atiyah $L^2$-index of Problem \ref{Integrality}   is well defined as 
$$
\Ind_{(2)} (\tD) \, \, := \, \, \dim_{(2)} \Ker (\tD^+) - \dim_{(2)} \Ker (\tD^-).  
$$
Using Theorem \ref{AlmostGap},  the following can be proved.

\begin{theorem}\cite{BenameurCovering}
Denoting by $\tau^{\av}$ and $\tau^{\reg}$ the average and regular traces as above, the following compatibility relations hold
$$
\tau^{\av}_* (\Ind (\maD_m)) = \Ind (D)\; \text{ and } \; \tau^{\reg}_* (\Ind (\maD_r)) = \Ind_{(2)} (\tD),
$$
\end{theorem}

Therefore, in the torsion-free case and if the index class $\Ind (\maD_r)$ belongs to the range of the BC map then $\Ind_{(2)} (\tD)$ is an integer and hence the pair $(M, D)$ is not a candidate for Problem \ref{Integrality}. This is due again  to the Atiyah theorem for closed manifolds. It is though suspected  that the  higher index class for non compact PSC near infinity complete spin manifolds can actually   escape the range of $\mu_\Gamma$. One probably needs to go beyond the cylindrical ends case of the APS index problem, which corresponds to $M=X\cup_{\pa X} (\pa X\times \R_+)$. This latter manifold  is obtained by the APS construction from a compact manifold with boundary $(X, \pa X)$ with a metric  whose restriction to the boundary  has PSC and which is product metric near the boundary, say by attaching a cylinder $\pa X\times \R_+$, and assuming  the compatibility of the $\Gamma$-coverings. In this  case, the Xie-Yu class  corresponds to a higher version of the  APS-index class \cite{APS1}. We may then ask the following question: \\

\begin{problem}\label{RangeMu1}
With $M$ being the APS cylindrical manifold  as above and $\Gamma$ torsion free, does the APS index class always belong to the range of the maximal BC assembly map.
\end{problem}

This question is of course about groups which do not satisfy the maximal BC conjecture. This special APS case is interesting because of its relation with reduced eta invariants and some other open questions that we shall discuss below.

In connection with the problems \ref{Integrality} and \ref{RangeMu1}, it is  unclear whether the Atiyah theorem, say the equality $ \Ind (D)= \Ind_{(2)} (\tD)$, is valid for general non-compact complete manifolds when one assumes only that $\Gamma$ is torsion free. But we expect it to be true in the APS case. This Atiyah-type theorem is  false when $\Gamma$ has torsion and is known to be true when $\Gamma$ is torsion-free and the class $\Ind (\maD_m)$ belongs to the range of the maximal BC assembly map, see  \cite{BenameurCovering}. Obviously, a positive answer to Problem \ref{Integrality} would show that the Atiyah theorem is false in general even if $\Gamma$ is torsion-free. Let us state though the following from \cite{BenameurCovering}:

\begin{problem}\label{Atiyah}
Find a pair $(M, D)$ with $\Gamma$ torsion free and the generalized Dirac operator invertible near infinity, such that $\Ind_{(2)} (\tD) - \Ind(D) \neq 0$. 
\end{problem}

\medskip

In order to find good candidates for the problems  \ref{Integrality} and \ref{Atiyah}, one needs to play with complicated ends but also with groups not satisfying the maximal BC conjecture.
When $\Gamma$ has torsion, we know already in the APS cylindrical case that $\Ind_{(2)} (\tD) - \Ind(D)$ coincides with the Cheeger-Gromov rho invariant and can be non zero. 
It would be interesting to investigate this defect invariant for the spin-Dirac operator under PSC near infinity slightly beyond the class of cylindrical ends manifolds. 
Recall that the Cheeger-Gromov invariant $\rho_{(2)} (D)=  \eta (D)- \eta_{(2)} (\tD)$ is the difference of the eta invariant of the spin-Dirac operator $D$ on $N$ and the $L^2$ eta invariant of its lift to the universal covering $\tN$ \cite{CheegerGromov}, with:
 $$
\eta_{(2)} (\tD) := \int_0^\infty \Tr_{(2)}(\tD e^{-t{\tD}^2})\,\frac{dt}{{\sqrt {\pi t}}} \; \text{ and }\;   \eta (D) :=  \int_0^\infty \Tr (D e^{-t{D}^2})\,\frac{dt}{{\sqrt {\pi t}}}.
 $$
There is an availlable APS $L^2$ formula  on $\tM$ \cite{Ramachandran}. 
Therefore, an affirmative answer to Problem \ref{RangeMu1} implies an affirmative answer to the following  problem:

\begin{problem}\label{RangeMu2}
Denote by $\rho_{(2)} (D)$ the Cheeger-Gromov  $L^2$ reduced eta invariant  of the spin-Dirac operator $D$ on a closed spin odd-dimensional manifold $N$ having a PSC metric. Assume that $\pi_1N$ is torsion free, then is it always true that $\rho_{(2)} (D)=0$?
\end{problem}

If the answer to Problem \ref{RangeMu2} is negative, then we expect again that a counter-example  will arise from  torsion-free groups which do not satisfy the maximal BC conjecture. 
Indeed, when the maximal Baum-Connes assembly map is an isomorphism, the answer to Problem \ref{RangeMu2} is affirmative, although this is not straightforward. See \cite{Keswani} and \cite{PiazzaSchick2005}.  

A fruitful   approach which significantly improved our understanding of the closed relation between the Cheeger-Gromov invariant and the maximal BC map, was introduced  by Higson and Roe in their series of papers \cite{HigsonRoeAnalysis1, HigsonRoeAnalysis2, HigsonRoeAnalysis3}. They constructed a six-term exact sequence encompassing the BC map with an explicit defect structure group, and they could deduce  in \cite{HigsonRoe2008} that for torsion-free groups,  the vanishing of their structure group implies the vanishing of the APS reduced eta invariants.  The same result is true for  the more general case of the Cheeger-Gromov $L^2$ rho invariant involved in Problem \ref{RangeMu2}, but this needed more techniques from semi-finite von-Neumann algebras, we recall the precise result below.  \\
 
In the Higson-Roe sequence, the structure group, denoted $\maS_*(\Gamma)$,   is an analytic version of the structure group appearing in the fundamental surgery exact sequence.  The Higson-Roe periodic six-term exact sequence,  can be summarized by the following triangle (we chose the maximal completions here):
 \begin{displaymath}\label{Figure1}
\xymatrixcolsep{1pc}\xymatrix{
K_{*}(B\Gamma) \ar[rr]^{\mu_\Gamma}  &  & K_{*}(C^*\Gamma)  \ar[dl]^{}\\ & \maS_{*+1} (\Gamma) \ar[ul]_{}
}\\
\end{displaymath} 
where $\mu_\Gamma$ is the BC assembly map and where 
 the groups $\maS_1(\Gamma)$ and $\maS_0(\Gamma)$ are defined analytically and vanish precisely when $\mu_\Gamma$ is an isomorphism.   There is a similar exact sequence for the regular completions which is easier to construct.

Let us explain  the relation of Problem \ref{RangeMu2} with the Higson-Roe work. So we assume that $\Gamma$ is torsion free. The above sequence is actually obtained by taking inductive limits in exact sequences for smooth Galois $\Gamma$-coverings  $\tM\to M$ given as follows:
 $$
\cdots \longrightarrow  K_0(M) \longrightarrow K_0(C^*\Gamma) \longrightarrow  \maS_1^\Gamma (\tM) \longrightarrow K_1(M) \longrightarrow K_1(C^*\Gamma) \longrightarrow  \cdots
 $$
In the odd dimensional spin and PSC case, the index class of the Dirac operator vanishes in $K_1(C^*\Gamma)$, therefore its class $[D]\in K_1(M)$ is the image of some class in $\maS_1^\Gamma (\tM)$. This latter, can explicitely be represented using the functional calculus of the operator $\maD_m$   \cite{HigsonRoe2008},  and its pushforward class in $\maS_{1} (\Gamma)$ is the higher rho class, and denoted $\rho (\maD_m)$. An obvious observation is that when the maximal BC map is an isomorphism, the higher rho class vanishes by construction. 
The difficult part in this approach is to be able to extract the numerical Cheeger-Gromov  invariant  from the rho class.   We have the following

\begin{theorem}\cite{BenameurRoyJFA}
Recall that $\Gamma$ is torsion free. There exists a group morphism $\xi_{(2)}: \maS_1(\Gamma)\rightarrow  \R$ which sends,  for any spin closed odd dimensional manifold $M$ with PSC, the higher rho class $\rho (\maD_m)$ on $M$, to the Cheeger-Gromov $L^2$ invariant, i.e.  $
\xi_{(2)} \left(\rho (\maD_m)\right) = \rho_{(2)} (\tD).$  
\end{theorem}

So this explains why when  the maximal BC map is an isomorphism, we necessarily have $\rho_{(2)} (\tD)=0$. Problem \ref{RangeMu2} is hence a consequence of the vanishing of the higher rho class in $\maS_1(\Gamma)$. The above theorem generalizes the similar one proved for the APS invariant in \cite{HigsonRoe2008}. 

\begin{remark}
When interested in homotopy invariance and the signature operator, the above theorem applies again. Namely there is a $\rho$-class  associated with any oriented homotopy equivalence between two oriented odd manifolds, whose image under  $\xi_{(2)}$  is the difference of the Cheeger-Gromov $L^2$ rho invariants of the  two oriented manifolds. This gives the oriented homotopy invariance of the Cheeger-Gromov $L^2$ invariant for the signature operator. 
 \end{remark}
 
\begin{remark}
When $\Gamma$ has torsion, the behavior of the Cheeger-Gromov rho invariant is totally different since it  is in general non trivial and it has been used to distinguish non-concordant classes of PSC metrics, see  \cite{BotvinikGilkey}. For the signature Cheeger-Gromov invariant, it is also no more homotopy invariant, and it was  used in \cite{ChangWeinberger} to prove that if  $\dim (M)=4k-1$ with $k\geq 2$ and $\pi_1M$ has torsion, then there exist infinitely many smooth manifolds which are homotopy equivalent to $M$ but not diffeomorphic to $M$.
\end{remark}

A possible approach to better understand Problem \ref{RangeMu2} stated above is to try to cook up a better geometric structure group, going beyond the  realm of $C^*$-algebras. 
%
%
%
The higher rho class is actually closely related with the Gromov-Lawson index in the cylindrical ends case, although there is no available APS higher index theorem so far. Still we can use the so-called delocalized APS theorem proved by Piazza and Schick to deduce this close relation for the two higher invariants. Indeed, for a Galois $\Gamma$-covering of spin even dimensional manifold with boundary $(\tM, \pa\tM)\to (M, \pa M)$, with PSC at the boundary,  one has the two classes recalled above. Namely, the higher index class $\Ind (\maD_M)\in K_0 (C_r^*\Gamma)$ and the higher  class  $\rho(D_{\pa M}) \in  \maS_1^(\Gamma (\pa\tM))$ whose pushforward in $\maS_1(\Gamma)$ is the $\rho (\maD_{\pa M})$. These two classes actually correspond to each other under functoriality. More precisely, the following compatibility theorem was  proved by Piazza and Schick

\begin{theorem}\cite{PiazzaSchick2014}
Denote by $i_*:\maS_1^\Gamma (\pa\tM)\rightarrow \maS_1^\Gamma (\tM)$ the natural map given by functoriality and by $j_*: K_0 (C^*\Gamma)\rightarrow \maS_1^\Gamma (\tM)$ the map in the Higson-Roe exact sequence for the $\Gamma$-covering $\tM\to M$. Then one has
$$
i_* \left(\rho(D_{\pa M}\right)  = j_*\left(\Ind (\maD_M)\right).
$$
\end{theorem}
 This result has been extended to all dimensions by Xie and Yu in \cite{XieYu2014b}.
In fact, this result was used in \cite{PiazzaSchick2014} and in \cite{XieYu2014b} to prove a stronger theorem providing a bridge between the Stolz exact sequence \cite{StolzSeq} and the Higson-Roe exact sequence \cite{HigsonRoeAnalysis3}.  We deduce that the vanishing of the rho class of the boundary operator on $\pa M$ implies the vanishing of the higher Cheeger-Gromov index class once pushed in the structure group, a coarse-type vanishing result. 

%
%
%
%

%

\begin{remark}\
There exists a Higson-Roe exact sequence for actions of discrete countable groups on (compact) spaces \cite{BenameurRoyII}. These actions correspond to a first class of foliations given by suspended flat bundles, for which the Cheeger-Gromov rho invariants and the higher indices were studied in \cite{BenameurPiazza}. We devote the next section to the construction of the higher index class of Xie-Yu for foliations.
\end{remark}

\section{Complete foliations and the $K$-theory index}

There is a similar BC assembly map for \'etale groupoids and especially for foliation groupoids. We shall devote this section to the setting of foliations. As explained for groups, higher index classes associated with spin-Dirac operators under the condition of PSC near infinity play a central part in the Gromov-Lawson approach through relative index theory, but they also provide exotic elements that may escape the range of the  BC map. 
We explain here how to assign  to any  foliated manifold with a leafwise spin structure  and with a complete metric whose leafwise scalar curvature is uniformly positive off some compact subspace of $M$, a spin-index class living in the $K$-theory of the ``space of leaves''. We involve in this definition any coefficient hermitian bundle whose leafwise curvature is small enough near infinity, in particular the $K$-theory index class is well defined for the leafwise Dirac operator with coefficients in any bundle whcih is leafwise flat off some compact subspace of $M$. Notice that the intersection of a compact subspace with a fixed leaf is not relatively compact in that leaf in general, so the restriction of the Dirac operator to a given leaf will not have a well defined Fredholm index. This is not surprising and the problem already appears for closed manifolds in the construction of the index by Connes \cite{ConnesIntegration}.
For the reader who is unfamiliar with groupoid foliations and their $C^*$-algebras, we have  rephrased in  the next section most of the constructions in the simpler case of locally trivial fibrations with a groupoid-free language. We refer the reader to the excellent reference \cite{Lawson} for the  geometric constructions with foliations, see also  \cite{Moerdijk}. 

There is a choice here for the $C^*$-algebra of the foliation, and we shall use a generating groupoid of the foliation which can equally be the holonomy groupoid or the monodromy groupoid. This latter giving more interesting invariants as we shall explain.  
We hence assume that $(M, F)$ is a smooth foliation with a complete metric $g$ which induces a complete metric on the even dimensional leaves and whose leafwise scalar curvature is uniformly positive off some compact subspace of $M$.  The leafwise tangent bundle $F$ is given at any $m\in M$ by the vector subspace of $T_m M$ composed of  vectors which are tangent to the leaf $L_m$ through $m$. The spin bundle associated with the fixed spin structure on  $F$ is denoted $S$ and it is equipped as usual with its hermitian structure (linear for the second entry and antilinear for the first) and spin connection, all compatible with $g$ and allow to define the leafwise Dirac operator $D$. The groupoid of $(M, F)$ is denoted $\maG$  and is abusively identified with its space of arrows. $\maG$ can be assumed to be Hausdorff for simplicity. The source and range maps of $\maG$ are denoted $s$ and $r$ respectively, and we shall use the common notation 
$$
\maG_A:=s^{-1} (A), \maG^B:= r^{-1} (B)\text{ and }\maG_A^B:=s^{-1} (A)\cap r^{-1} (B),\quad \text{ for }A,B\subset M.
$$
Notice that $\maG_A^A$ is a subgroupoid of $\maG$. Given $m\in M$, the holonomy/monodromy covering of the leaf $L_m$ through $m$ will be the $\maG_m^m$-covering $r:\maG_m=\maG_{\{m\}} \to L_m$. The operator $D$ being a leafwise differential operator it restricts to the leaves. Morever, we shall rather consider for any $m$, its lifts $D_m$, a Dirac operator on the manifold $\maG_m$ acting on the $\maG_m^m$-equivariant spin bundle $r^*S$. The family $(D_m)_{m\in M}$ is then a globally $\maG_T^T$-equivariant family of Dirac operators, see \cite{ConnesIntegration}.  Let us choose a locally finite distinguished good open cover for the foliation $(U_i)_{i\in I}$ by relatively compact open submanifolds  with small transversal $T_i\in U_i$ so that $U_i\simeq T_i\times \D^p$ with $\D^p$ the unit disc in $\R^p$. We also assume as usual  that $\bar T_i\cap \bar T_j=\emptyset$ for $i\neq j$. The union $\cup_i T_i$ is then a complete transversal submanifold, meaning that it is a transversal submanifold which cuts all the leaves.  Using the complete transversal $T$, we may replace $\maG$ by the reduced \'etale groupoid $\maG_T^T$ without changing the Morita equivalence class of $\maG$ \cite{HilsumSkandalis}, see also \cite{BenameurPacific}.  The reduced $C^*$-algebra associated with  the \'etale grooupoid $\maG_T^T$ is denoted $C^*\maG_T^T$, it is the completion of the involutive convolution algebra $C_c^\infty (\maG_T^T)$ with respect to its regular representation in the continuous field of Hilbert spaces  $(\ell^2\maG^T_m)_{m\in T}$, see \cite{Renault}. The groupoid $\maG_T^T$ acts freely and properly on the manifold $\maG_T$  with quotient diffeomorphic to the ambiant manifold $M$, it actually embodies the Morita equivalence between the two groupoids $\maG$ and $\maG_T^T$, see \cite{HilsumSkandalis} and also \cite{BenameurPacific}. 

The space $C^\infty_c (\maG_T, r^* S)$  of   smooth compactly supported sections of $r^*S$ over $\maG_T$ is then a right prehilbertian $C^\infty_c(\maG_T^T)$-module with a left representation of the algebra $C_c^\infty (\maG)$. Using the Haar system $m\mapsto \nu_m$ on $\maG_m$ lifted from  a Lebesgue class measure on the leaves, the rules are
$$
(\xi f) (\gamma) := \sum_{\beta\in \maG_{s(\gamma)}^T} \xi (\gamma\beta^{-1} f(\beta)\; \in S_{r(\gamma)}\; \text{ and }\; \langle \xi, \eta\rangle (\beta):=\int_{\gamma_1\in \maG_{r(\beta)}} \langle \xi (\gamma_1), \eta (\gamma_1\beta)\rangle_{S_{r(\gamma_1)}} d\nu_{r(\beta)} (\gamma_1),
$$
for $\xi\in C_c^\infty (\maG_T, r^* S)$, $f\in C_c^\infty (\maG_T^T)$, $\beta\in \maG_T^T$ and $\gamma\in \maG_T$. 

Let us  denote by $\maE$ the completion of $C^\infty_c (\maG_T, r^* S)$ with respect to this inner product \cite{BenameurRoyPoincare}.  Then $\maE$ is a Hilbert module over the $C^*$-algebra $C^*\maG_T^T$ and the representation of $C_c^\infty (\maG)$  extends to the representation $\pi: C^*\maG \to \End_{C^*\maG_T^T} (\maE)$  which is valued in the $C^*$-algebra $\maK (\maE)$ of compact operators \cite{KasparovStinespring}.

The leafwise spin connection $\nabla$ can be used to define the   Sobolev Hilbert modules $\maE_k$ for any $k\geq 0$. In terms of continuous fields of Hilbert spaces in the sense of Dixmier, $\maE_k$ corresponds to the field of Sobolev spaces of $\maG_m$ with coefficients in $r^*S$. Said differently, $\maE_k$ is obtained by  completing  the space $C_c^\infty (\maG_T, r^*S)$ with respect to the norm of $C^*\maG_T^T$ and to the $C_c^\infty (\maG_T^T)$-valued  inner product
$$
\langle \xi_1, \xi_2\rangle_k:= \sum_{j=0}^k \langle \nabla^j \xi_1, \nabla^j \xi_2\rangle.
$$
The inclusion $\maE_1\hookrightarrow \maE$ is then an adjointable operator.

If $E$ is an extra hermitian bundle over $M$ with a leafwise connection $\nabla^E$, then one may replace $\Dir$ by the $\nabla^E$-twisted Dirac operator $D$ and the Lichnerowicz formula for $D$ has the form \cite{LawsonMichelsohn}:
$$
D^2 \,\, = \,\, \nabla^* \nabla  +   \frac{\kappa}{4}  +  \mathcal{R}^E,
$$
with $\kappa:M\to \R$ being the scalar curvature of the restriction of the metric to the leaves, a smooth function, and  $\maR^E$ a leafwise curvature contribution of the connection $\nabla^E$. All the previous considerations then work again replacing $S$ by $S\otimes E$. We shall denote again by $\maE$ the same Hilbert module but for this latter twisted spinor bundle. 

Assume that there exists a constant $\kappa_0 >0$ such that  $\frac{\kappa}{4} - \vert \maR^E\vert \geq \kappa_0>0$   outside  some compact subspace of $M$. This happens for instance when $\kappa$ is uniformly positive off some compact subspace and when $E$ is leafwise flat. Leafwise flat bundles were sometimes called  foliated bundles in some references, they have basic characteristic classes.  An even more interesting situation is when $E$ is only flat along the leaves outside a compact subspace, or even only an almost flat $K$-theory class. 
Then  we fix   a compactly supported non-negative  smooth function  $\rho$ on $M$ such that $\kappa- 4\vert \maR^E\vert  + 4\rho \geq 4\kappa_0 >0$ over $M$. The  operator $D^2+\rho$ is then regular and has bounded (self-adjoint) inverse $(D^2+\rho)^{-1}$, whose  square root $(D^2+\rho)^{-1/2}$ is defined using the spectral theorem for Hilbert modules. 
The following lemma was stated in \cite{XieYu2014} for galois coverings and the proof is almost the same.

\begin{lemma}\
For any  smooth  compactly supported function $f:M\to \C$, the operator $f(D^2+\rho)^{-1/2}$ is a compact operator on the Hilbert module $\maE$. 
\end{lemma}

\begin{proof} 
 The  operator $(D^2+\rho)^{1/2}$ is a regular operator which is injective with bounded inverse, this can be justified exactly as   in  \cite{XieYu2014} for the Galois covering case, see also \cite{BenameurCovering}. More precisely,  let $\sigma\in \Dom (D^2)\subset \maE$. We only give the proof for trivial one dimensional $E$ and ignore $\maR^E$ for simplicity. Then 
$$
\langle (D^2+\rho)\sigma, \sigma\rangle = \langle (\nabla^*\nabla  + (\frac{\kappa}{4} +\rho) \Id) \sigma, \sigma\rangle\geq  \vert\vert \nabla \sigma \vert\vert^2 + \kappa_0 \vert\vert \sigma\vert\vert^2 \geq \kappa_0 \vert\vert \sigma\vert\vert^2.
$$
Hence the operator $D^2+\rho$ is injective and has  a bounded inverse. The spectral theorem on Hilbert modules again identifies the inverse with $(D^2+\rho)^{-1/2}$. Indeed,  the spectrum of $(D^2+\rho)^{1/2}$ is contained in $[\sqrt{\kappa_0}, +\infty)$. Therefore, the operator $(D^2+\rho)^{-1/2}$ is an adjointable  non-negative operator on $\maE$. Moreover, for $\sigma\in \maE$ we can write
\begin{eqnarray*}
\vert\vert (D^2+\rho)^{-1/2} \sigma\vert\vert_1^2  & = & \langle (D^2+\rho)^{-1} \sigma, \sigma \rangle  + \langle (D^2+\rho)^{-1/2} \nabla^*\nabla (D^2+\rho)^{-1/2} \sigma, \sigma\rangle \\
& \leq  & \frac{1}{\kappa_0}  \vert\vert\sigma\vert\vert^2 + \langle (D^2+\rho)^{-1/2} [(D^2+\rho)-(\rho+\frac{\kappa}{4})] (D^2+\rho)^{-1/2} \sigma, \sigma\rangle
 \end{eqnarray*}
Since $D^2+\rho)^{-1/2} (D^2+\rho) (D^2+\rho)^{-1/2} = \Id$ and  $\rho +\frac{\kappa}{4} \geq \kappa_0$, we deduce
$$
\vert\vert (D^2+\rho)^{-1/2} \sigma\vert\vert_1^2  \leq (1+\frac{1}{\kappa_0}) \vert\vert\sigma\vert\vert^2.
$$
Therefore, the operator $(D^2+\rho)^{-1/2}$ is a bounded operator from $\maE$ to $\maE_1$ as announced. The Rellich lemma allows to conclude.
More precisely, if $f$ is a smooth compactly supported function on $M$, then pointwise multiplication by $f$  extends to a compact operator from $\maE_{1}$ to $\maE$. This can be proved by reducing to a local chart and by using \cite{FomenkoMiscenko}[Lemma 3.3].
%
\end{proof}

We call  Fredholm operators all  adjointable operators which are invertible module compact operators. 
We can now state the main result of this section.

\begin{theorem}\label{Fredholm}
Assume that the  metric $g$ satisfies $\frac{\kappa}{4} - \vert\maR^E\vert  \geq \kappa_0>0$ off some compact subspace $K$ of $M$ where $\kappa$ is the scalar curvature along the leaves as above. Then the odd adjointable operator $F_\rho:=(D^2+\rho)^{-1/2}D$ is an odd symmetry modulo compact operators, i.e.
\begin{itemize}
\item $F_\rho$ is odd for the grading;
\item $F_\rho^2 - \Id$ is a compact operator;
\item $F_\rho^*-F_\rho$ is a compact operator.
\end{itemize}
In particular, the chiral part $F_\rho^+$ has a well defined index class in $K_0(\maK(\maE))\simeq K_0(C^*\maG)$.
\end{theorem}

\begin{proof}
We follow \cite{XieYu2014} which can be adapted here. The spectral theorem shows that the operator $(D^2+\rho)^{-1/2}$ can also be defined using the  convergent integral in the $C^*$-algebra of adjointable operators on $\maE$:
$$
(D^2+\rho)^{-1/2} := \frac{2}{\pi} \int_0^{+\infty}\;  (D^2+\rho+\mu^2)^{-1}\;  d\mu.
$$
Hence 
$$
[D, (D^2+\rho)^{-1/2}]=\frac{2}{\pi} \int_0^{+\infty} (D^2+\rho+\mu^2)^{-1} [\rho, D] (D^2+\rho+\mu^2)^{-1} \, d\mu,
$$
and the operator $[\rho, D]$ being given by Clifford mutiplication with the leafwise gradient of a compactly supported  smooth function, it is an adjointable operator. Since  $D^2+\rho\geq \kappa_0 \Id$, we deduce
$$
\vert\vert (D^2+\rho+\mu^2)^{-1} [\rho, D] (D^2+\rho+\mu^2)^{-1}\vert\vert \leq \vert\vert [\rho, D]\vert\vert \frac{1}{(\kappa_0+\mu^2)^2}.
$$
%
%
%
Moreover, $\rho$ and $[D, \rho]$ being both compactly supported zero-th order  operators, the Rellich lemma implies that the operator 
$$
[\rho, D] (D^2+\rho+\lambda^2)^{-1}
$$ 
is a compact operator. Therefore the commutator $[D, (D^2+\rho)^{-1/2}]$ is compact and so $F_\rho^*-F_\rho\in \maK (\maE)$.
Now,
\begin{eqnarray*}
F_\rho^2 &=& F_\rho (D^2+\rho)^{-1/2} D \\
& = &(D^2+\rho)^{-1/2} D^2 (D^2+\rho)^{-1/2} - F_\rho [D, (D^2+\rho)^{-1/2}] \\
&=& I - (D^2+\rho)^{-1/2} \rho (D^2+\rho)^{-1/2} - F_\rho [D, (D^2+\rho)^{-1/2}]
\end{eqnarray*}
Now, the operator $[D, (D^2+\rho)^{-1/2}] F_\rho$ is compact, and a similar verification shows  that the operator $(D^2+\rho)^{-1/2} \rho (D^2+\rho)^{-1/2}$ is also compact. Hence $F_\rho$ is Fredholm with parametrix given by $F_\rho$ itself. 

The above properties of $F_\rho$  allow by the standard boundary construction in $K$-theory, to assign to the operator $F_\rho$ an index class $\Ind (F^+_\rho)$ which lives in the $K$-theory group of the $C^*$-algebra $\maK (\maE)$. To be specific, the operator 
$$
F_\rho^+:= (D^+D^- +\rho)^{-1/2} D^+: \maE^+\longrightarrow \maE^-,
$$
is Fredholm with quasi inverse given by 
$$
F_\rho^-:=  (D^-D^+ +\rho)^{-1/2}D^-: \maE^-\longrightarrow \maE^+.
$$
If we denote by
$$
\maS^+ = \Id_{\maE^+} - F_\rho^-F_\rho^+\text{ and }\maS^- = \Id_{\maE^-} - F_\rho^+ F_\rho^-,
$$
then the  index class  can be represented by the difference $K$-class $[e]-[f]$,  where $e$ and $f$ are the two idempotents  given by:
$$
e:=\left(\begin{array}{cc} (\maS^+)^2 & \maS^+ F_\rho^- \\ F_\rho^+\maS^+ (\Id_{\maE^+}+\maS^+) & \Id_{\maE^-} -(\maS^-)^2 \end{array}\right) \text{ and } f=\left(\begin{array}{cc} 0  & 0 \\ 0 & \Id_{\maE^-} \end{array}\right)
$$
\end{proof}

The index class constructed in the previous theorem does not depend on the choice of the function $\rho$ satisfying $\kappa-\vert\maR^E\vert +4 \rho\geq 4\kappa_0 >0$ since two such choices lead to an operator homotopy between the cycles. Moreover, the Morita equivalence between $\maK (\maE)$ and the $C^*$-algebra $C^*\maG$ provides the following definition of the index class of our leafwise Dirac operator.

\begin{definition}
Assume as above that there exists $\kappa_0>0$ such that $\frac{\kappa}{4} - \vert \maR^E\vert \geq \kappa_0>0$ off some compact subspace of $M$, then the $E$-twisted leafwise Dirac operator $D$ has a well defined index class $\Ind (D)$  in the $K$-theory group of $C^*\maG$. More precisely, under the Morita  isomorphism $K_0(\maK (\maE))\simeq K(C^*\maG)$, the index class $\Ind (D)\in K(C^*\maG)$ is the class  corresponding to  $\Ind (F_\rho ^+)\in K_0(\maK(\maE))$. 
\end{definition}

When the manifold $M$ is compact and we use the Holonomy groupoid,  our index class coincides with the Connes-Skandalis index class. When the foliation is top-dimensional, we recover the Gromov-Lawson  integer index as defined in  \cite{GromovLawson3} when $\maG$ is the holonomy groupoid, and we recover the higher Gromov-Lawson index class in $K_0(C^*\pi_1M)$ as defined in \cite{XieYu2014} when $\maG$ is  the monodromy groupoid. 

Again and as in the case of a single complete manifold, there is a well defined Baum-Connes assembly map for the monodromy groupoid which is conjectured  to be an isomorphism and which is related with other conjecture such as the Novikov conjecture for foliations but also with a foliated version of the Gromov-Lawson-Rosenberg conjecture, and whose receptacle is precisely the $K$-theory group of the $C^*$-algebra of the monodromy groupoid of the given foliation.  As in Section \ref{SectionCovering} and even in the torsion-free case, it is not clear whether the index class constructed here belongs  to  the range of this assembly map.

In \cite{BenameurHeitsch21}, the authors constructed by a different method,  a Gromov-Lawson measured index in full generality when the foliation admits a holonomy invariant measure $\Lambda$, now using the kernel bundle of the leafwise Dirac operator. They also proved a measured relative index theorem generalizing the Connes measured index theorem \cite{ConnesIntegration} exactly as the Gromov-Lawson index theorem generalized the Atiyah-Singer index theorem. The relation with our index class is obtained using the induced group morphism $
\tau_\Lambda: K(C_r^*\maG)\longrightarrow \R.$

A particularly interesting situation is when the foliation is given by a suspension, as studied in \cite{BenameurPiazza}. Assume that $\Gamma$ is the fundamental group of a smooth even complete spin manifold $W$ and that $\Gamma$ acts by diffeomorphisms on a closed manifold $X$. Then the smooth manifold $M$ which is the quotient of $\tW\times X$ under the diagonal action $(\tw, x)\gamma:= (\tw\gamma, \gamma^{-1}x)$ is a bundle over $W$ with fibers diffeomorphic to $X$ and it has a smooth foliation which is transverse to the fibers with leaves being the projections in $M$ of the submanifolds $\tW\times\{x_0\}$ for $x_0\in X$. In this case one may choose a complete transversal as a specific fiber $X_w$ and restricting the monodromy groupoid to this transversal, one gets a groupoid which is isomorphic to the crossed product groupoid $X\rtimes \Gamma$, say the  groupoid associated with the action of $\Gamma$ on $X$. The above construction gives us in this case and under the above assumptions a class $\Ind (\maD_m)$  in  the $K$-theory group $K(C(X)\rtimes \Gamma)$ where we denoted by $C(X)\rtimes \Gamma$ the maximal $C^*$-algebra of this groupoid $X\rtimes \Gamma$. If now, $\mu$ is a $\Gamma$-invariant measure on $X$, then the average and regular traces associated with $\mu$  were defined in \cite{BenameurPiazza}, extending the case of groups described in the previous section. These traces induce the group morphisms
$$
\tau_\mu^{reg} : K(C(X)\rtimes \Gamma) \longrightarrow \R\text{ and } \tau_\mu^{av} : K(C(X)\rtimes \Gamma) \longrightarrow \R.
$$
We obtain by applying these morphisms, two measured indices $\Ind_\mu^{reg} (D)$ and $\Ind_\mu^{av} (D)$ for any $\Gamma$-equivariant (smooth) family over $X$ of metrics on $\tW$ which induces PSC near infinity in the quotient $M$, say such that the family of scalar curvature functions, a global $\Gamma$-invariant smooth function on $\tW\times X$ is a uniformly positive function  outside some  $\Gamma$-compact subspace of $\tW\times X$. The Atiyah theorem holds in this case, when the manifold $W$ is compact, see \cite{BenameurPiazza}. However and as observed in the previous section for the case $X=\{\bullet\}$, there is no evidence that such equality holds in general, even when $\Gamma$ is torsion-free. 

The measured rho invariant $\rho_\mu (X, \Gamma)$ was also defined and studied in the similar situation of $\Gamma$-actions in \cite{BenameurPiazza} with $\Gamma$ now the fundamental group of an odd closed manifold. It was proved there that when $\Gamma$ is torsion-free and the maximal BC map for $X\rtimes \Gamma$ is an isomorphism, this measured $\rho$ invariant vanishes again. Finally, the Higson-Roe exact sequence relating the BC assembly map with a structure group does also exist in this more general setting, as recently proved in \cite{BenameurRoyII}.

When the foliated manifold $(M, F)$ is obtained by attaching a foliated cylinder in an Atiyah-Patodi-Singer problem for a foliated compact manifold with boundary transverse to the leaves, with a metric which has Positive leafwise Scalar Curvature at the boundary and with the bundle $E$ having leafwise flat connection at the boundary, the above construction of the index class corresponds to the expected APS index class. Moreover, for foliations given by suspensions as studied in \cite{BenameurPiazza} and where $\maG_T^T$ becomes the crossed product groupoid $T\rtimes \Gamma$ corresponding to an action of $\Gamma$ on the compact manifold $T$, the corresponding Higson-Roe  exact sequence  holds as proved in \cite{BenameurRoyII}, with the corresponding structure group $\maS_1 (T\rtimes \Gamma)$ where the higher foliated rho class lives. This higher rho class is clearly related, for foliations with boundaries having leafwise PSC metrics at the boundary, with the index class constructed above.

 In the next section, we   give further developpements  in the locally trivial fibration case, where the Bismut-Cheeger higher formula is available. 

\section{Some GL invariants for families}

Let $\pi:M\to B$ be a fibre bundle with a smooth non-compact typical fiber, over a compact  space $B$. In the sequel, we shall always assume for simplicity that the fibers are connected and that the restriction of the projection to any complement of a compact subspace of the total space is surjective. So, $\pi$ is a submersion which is locally trivial with the local picture given by $U\times F$ for open subspaces $U$ and with $F$ being the typical fiber, a smooth manifold. The changes of trivialization will be assumed to be induced by homeomorphisms which restrict to smooth diffeomorphisms of the fiber $F$. We assume that the fibre bundle can be induced with  continuous in the $B$ variable, fiberwise smooth uniformly  (in the $B$-variable) complete metric along the fibers $g= (g_b)_{b\in B}$ as well as a fiberwise spin structure $\sigma=(\sigma_b)_{b\in B}$. The associated fiberwise spin-Dirac  operator is denoted $D=(D_b)_{b\in B}$ and it acts on the  spinor bundle $S$ of the fiberwise tangent bundle. The bundle $S$ is endowed with the usual  hermitian structure inducing the spin connection, see for instance \cite{LawsonMichelsohn}. We may and shall as usual assume  that the fibers are even dimensional so that $S=S^+\oplus S^-$ is $\Z_2$-graded. When the metric has fiberwise Positive Scalar Curvature off some compact subspace of the ambiant space $M$, then on any fiber $M_b$ the Gromov-Lawson theorem tells us that the Dirac operator $D_b$ has finite dimensional kernel. Hence, we find ourselves exactly as in the case of closed fibre bundles as studied by Atiyah and Singer in \cite{AtiyahSinger4}, and by the same Atiyah-Singer construction, the family $D=(D_b)_{b\in B}$ can be shown to have a well defined higher index class in the  $K$-theory group of $B$. This class, up to the Atiyah-Singer   perturbation of $D$ \cite{AtiyahSinger4}, is  the class of the so-called index bundle.  

Let us review the construction of the index class from the previous section in this simple foliation cas and as explained above,  only assuming transverse continuity. Denote by $\Gamma_c (M, S)$ the space  of compactly supported continuous and fiberwise smooth sections of $S$. Let then $\maE$ be the $\Z_2$-graded Hilbert $C(B)$-module which is the completion of $\Gamma_c (M, S)$ with respect to the fiberwise inner product, i.e.
$$
\langle \xi_1, \xi_2\rangle (b):= \int_{M_b} \langle \xi_1 (m), \xi_2(m)\rangle_{S_m} dm
$$
where $dm$ is the induced  fiberwise Lebesgue class measure. The fiberwise spin Dirac operator  is denoted $\Dir:\Gamma_c (M, S)\to \Gamma_c (M, S)$, and we may again twist it by a hermitian bundle $E$ with a fiberwise connection $\nabla^E$ to defined the fiberwise Dirac operator $D$. By classical arguments, $D$ extends to a regular self-adjoint, odd for the $\Z_2$-grading, operator on $\maE$, and $D^2$ also extends to a non-negative regular operator on $\maE$.

Again, we assume that the fiberwise metric has uniformly positive fiberwise scalar curvature off some compact subspace of $M$ and moreover  that the curvature of $\nabla^E$ is small enough off the same compact subspace, so that as in the previous section we may  choose   a compactly supported non-negative continuous and fiberwise smooth function  $\rho$ on $M$ such that the  operator $D^2+\rho$ has bounded self-adjoint inverse $(D^2+\rho)^{-1}$. The square root of this latter operator is the operator $(D^2+\rho)^{-1/2}$. Notice that the conditions imposed on $E$ are for instance satisfied for any fiberwise flat bundle, or even when it is only fiberwise flat outside some compact subspace.

Theorem \ref{Fredholm} then becomes here:

\begin{theorem}\label{Fredholm}
Under the above assumptions, the operator  $F_\rho:=(D^2+\rho)^{-1/2}D$ is an odd symmetry modulo the ideal of compact operators of $\maE$, and the index of its chiral part is  a well defined class in $K_0(\maK(\maE))\simeq K(B)$.
\end{theorem}

The  index class  can again be represented by the difference $K$-class $[e]-[f]$ as in the previous section. 

\begin{definition}
Under the above assumptions on the fiberwise metric and the fiberwise connection of $E$,   the fiberwise $E$-twisted Dirac operator $D$ has a well defined index class $\Ind (D)$  in the $K$-theory group of the base space $B$. 
\end{definition}

When the index bundle of $D$ exists, one can prove that the index class is given by this index bundle, i.e.
$$
[\Ker (D^+)]-[\Ker (D^-)]\;  \in \; K(B).
$$ 
As a first application, let us relate our higher index class $\Ind (D)$ with the family APS-index class as studied in \cite{BismutCheeger1, BismutCheeger2, MelrosePiazza1, MelrosePiazza2}. Assume then that $M\to B$ is a smooth fibre bundle by smooth even spin manifolds with boundary over a smooth manifold $B$, so that $\pa M \to B$ is also a smooth fibre bundle. We fix a hermitian bundle $E$ over $M$ with a hermitian connection $\nabla^E$ whose restriction to the boundary induces a connection with fiberwise small curvature, for instance fiberwise flat at the boundary. As in  \cite{BismutCheeger1}[pp 341-342], we assume that the fiberwise metric is of product type $g+dx^2$ in a vertical collar neighborhood of $\pa M$ in $M$, where $g$ is the induced fiberwise metric on $\pa M\to B$ and has PSC, and the same assumptions are taken for the bundle $E$, so that the boundary fiberwise $E$-twisted Dirac operator $D_{\pa M}$ is invertible. In this case the Bismut-Cheeger index formula holds for the Chern character of the index of the family of Dirac operators with the fiberwise APS boundary conditions involving in addition to the usual Atiyah-Singer integral a contribution from the boundary; the so-called eta form, see again \cite{BismutCheeger2}. It is worth ppointing out that even when the operator is not invertible at the boundary but has a trivial class in $K^1(B)$, one still has a higher APS formula which uses spectral sections, see \cite{MelrosePiazza1}. 

By the results of \cite{MelrosePiazza1},  the Bismut-Cheeger higher APS-index class $\Ind (D_M, APS)$  identifies  then with our higher index class as obtained  for the extended family of Dirac operators on the augmented fibre bundle $\what{M}\to B$ obtained by gluing a cylinder $\pa M\times [0, +\infty)\to B$ to the fibre bundle $M\to B$. Therefore, one may use the Bismut-Cheeger APS-index formula for its Chern character that we denote by $\iota (M\vert B, g, \nabla^E)$, which belongs to $H^{[0]} (B, \R)$, to deduce  its independence of  the extended metric and the extended connection of $E$.

In the same way, given a closed fibre bundle by odd dimensional spin manifolds $N\to B$ over a manifold $B$, with two metrics $g_0$ and $g_1$ which both yield complete fiberwise metrics with PSC, we may endow the complete fibre bundle $M=N \times \R\to B$ with any metric which coincides with $g_0$ on $N\times (-\infty, 0]$ and with $g_1$ on $N\times [1, +\infty)$. Adding the extra hermitian bundle $E$ over $N$ with two fiberwise flat connections $\nabla^0$ and $\nabla^1$,  we construct the resulting fiberwise Dirac operator $D_M$ on $M\to B$ which then has a well defined higher index class in $K(B)$. This is the exact generalization of the Gromov-Lawson invariant for the pair of PSC metrics $(g_0, g_1)$ on $N$, see  \cite{GromovLawson3}. We denote its Chern character in the even real cohomology of $B$ by 
$$
\iota \left(N\vert B; (g_0, \nabla^0);  (g_1, \nabla^1)\right)\; \in \; H^{[0]}(B, \R)\quad \text{ and without the bundle $E$: }\iota (N\vert B; g_0, g_1).
$$

\begin{remark}
The invariant $\iota \left(N\vert B; (g_0, \nabla^0);  (g_1, \nabla^1)\right)$ belongs to $H^{[0]}(B, \Q)$ but only  its image in $H^{[0]}(B, \R)$ can be computed by the Bismut-Cheeger formula.
\end{remark}

An obvious observation is that the invariant  $\iota \left(N\vert B; (g_0, \nabla^0);  (g_1, \nabla^1)\right)$ vanishes whenever there exists a metric of fiberwise PSC on $N\times [0, 1]$  and a fiberwise flat connection on $E\times [0, 1]$ which extend the given data at the boundary. So in particular, $\iota (N\vert B, g_0, \nabla^0; g_1, \nabla^1)$ vanishes when $g_0$ and $g_1$ belong to the same connected component of the  space of metrics of fiberwise PSC on $N\to B$ and $\nabla^0, \nabla^1$ are connected by fiberwise flat connections. Again by  \cite{MelrosePiazza1}, the invariant $\iota \left(N\vert B; (g_0, \nabla^0);  (g_1, \nabla^1)\right)$ coincides with  the Bismut-Cheeger higher APS-index for the fibre bundle with boundary $N\times [0, 1]\to B$ with the metric on the boundary which is $g_0$ on $N\times \{0\}$ and $g_1$ on $N\times\{1\}$ and any connection on $E\times [0, 1]\to N\times [0, 1]$  extending $\nabla^0$ and $\nabla^1$.  Therefore, the Bismut-Cheeger formula holds again and allows to prove the Gromov-Lawson cocycle relation

\begin{proposition}\ One has $
\iota \left(N\vert B; (g_0, \nabla^0);  (g_1, \nabla^1)\right) + \iota \left(N\vert B; (g_1, \nabla^1); (g_2, \nabla^2)\right)=\iota \left(N\vert B; (g_0, \nabla^0);  (g_2, \nabla^2)\right).$
\end{proposition}
One expects  that this cocycle relation already holds for the Gromov-Lawson higher index in $K(B)$, before taking the Chern character, but this needs a proof \`a la cut-and-paste like the one given in \cite{GromovLawson3}.

Following Gromov and Lawson, we  then consider for $\pi_N:N\to B$ as above the group $\Diff^{\infty, 0} (N\vert B)$ of homeomorphisms of $N$ which commute with the projection $\pi_N$ and induce  smooth diffeomorphisms of the fibers. If we  choose a fiberwise metric $g$ with Positive Scalar Curvature as above,  any $F\in \Diff^{\infty, 0} (N\vert B)$ gives a new fiberwise PSC metric $F^*g$. Denoting by $\Gamma (N\vert B)$ the component group of $\Diff^{\infty, 0} (N\vert B)$ we have the exact generalization of Theorem 4.48 and Corollary 4.49 in \cite{GromovLawson3}, the proof is exactly the same.

\begin{proposition}
The map $F\mapsto \iota (N\vert B; g,  F^*g)$ induces a well defined group morphism 
$$
\iota_g : \Gamma (N\vert B)\longrightarrow H^{[0]} (B, \R).
$$
\end{proposition}


\medskip
We end this section with a few remarks related with a family version of the Kreck-Stolz invariant constructed in \cite{KreckStolz}. We then need the following Green-Stockes lemma:

\begin{lemma}\label{Stockes}
Assume that $M\to B$ is a compact fibre bundle with boundary, and such that $\pa M\to B$ is  a subbundle. Assume that $\alpha$ and $\beta$ are closed differential forms on $M$ with restrictions to the boundary wchich are exact. Choose a differential form $\what\alpha$  on $\pa M$ such that  $\alpha\vert_{\pa M} = d\what\alpha$. Then
$$
\int_{M\vert B} \alpha\wedge \beta - \int_{\pa M\vert B} \what\alpha\wedge \beta\vert_{\pa M}.
$$
is a closed differential form on $B$. Moreover, if we denote  by $\tilde\alpha$ and $\tilde\beta$ any closed forms on $M$ which are cohomologous to $\alpha$ and $\beta$ respectively and vanish at the boundary, then  the following formula holds  {\underline{modulo exact forms}} 
$$
\int_{M\vert B} \alpha\wedge \beta - \int_{\pa M\vert B} \what\alpha\wedge \beta\vert_{\pa M} = \int_{M\vert B} \tilde\alpha\wedge\tilde\beta.
$$
\end{lemma}

\begin{proof}
We only need to show the last equality modulo exact forms since the RHS is a closed differential form. Let us choose $\alpha_0$ and $\beta_0$ in $\Omega^* (M)$ such that 
$$
\tilde\alpha-\alpha= d\alpha_0\text{ and } \tilde\beta-\beta= d\beta_0.
$$
The Green-Stockes theorem for fibrations with boundary which reads for any $\delta\in \Omega^* (M)$:
$$
\int_{M\vert B} d\delta - d\int_{M\vert B} \delta = \pm \int_{\pa M\vert B} \delta\vert_{\pa M}.
$$ 
This can be  proved locally and hence for a trivial fibration. The proof of the lemma is then  an easy direct computation starting by replacing $\alpha$ and $\beta$ by $\tilde\alpha-d\alpha_0$ and $\tilde\beta-d\beta_0$, applying the Green-Stockes formula. One has to be careful about the sign appearing in the Green-Stockes formula. 
\end{proof}

 Recall,  for a given spin fibre bundle $N\to B$ with  a $4k-1$-dimensional fibre and ambiant metric $g$ on $N$, the Bismut-Cheeger eta form $\what{\eta} (D_N)\in \Omega^{[0]}(B)$  associated with the Bismut superconnection derived from the  Levi-Civita connection on $(N, g)$ and the associated fiberwise spin connection, as in \cite{BismutCheeger1}. Assume  the  metric $g$ has Positive fiberwise Scalar Curvature and $E$ is fiberwise flat.  Then the following transgression formula  was proved in  \cite{BismutCheeger1, BismutCheeger2}:
$$
d\what{\eta} (D_N) = \int_{N\vert B} \what{A} (TN\vert B, g) \Ch (\nabla^E).
$$
In fact, there is a similar formula when one only assumes that the family  $D_N$ has trivial index class in $K^1(B)$ \cite{MelrosePiazza1}. In particular, this relation then holds for the untwisted spin-Dirac operator $\Dir_N$ under PSC along the fibers, and we have $d\what{\eta} (\Dir_N) = \int_{N\vert B} \what{A} (TN\vert B, g)$ in this case. 

Moreover, when the higher odd signature operator has trivial $K^1$-class, there exists a similar transgression formula for the signature operator $D_N^{sgn}$ which takes the following form:
$$
2 \, d\what{\eta} (D_N^{sgn}) = \int_{N\vert B} \LL (TN\vert B, g).
$$
An even differential form satisfying this relation is not unique of course and one can add any closed even form. Let us assume that the differential form $\what\eta (D_N^{sgn})$ is precisely  the eta form appearing the Bismut-Cheeger APS index formula, exactly as for the Dirac operator, this is true under the appropriate assumptions, see \cite{BismutCheeger2}.

We shall use the following notations from \cite{KreckStolz}:
\begin{itemize}
\item If $\alpha$ and $\beta$ are exact forms on $N$, then denote abusively $\int_{N\vert B} d^{-1}(\alpha\wedge \beta) := \int_{N\vert B} \what\alpha\wedge \beta$, for any differential form $\what\alpha$ such that $\alpha=d\what\alpha$. See Lemma \ref{Stockes}.
\item Set $a_k:=\frac{1}{2^{2k+1}(2^{2k-1}-1)}$.
\end{itemize}

We are now in a position to  introduce the higher Kreck-Stolz invariant $\sigma (N\vert B, g)$ in  the even real cohomology space $H^{[0]} (B, \R)$.

\begin{definition}
Assume that  the  fibre bundle $N\to B$  is spin with a $4k-1$-dimensional typical fiber and that the metric has Positive  fiberwise Scalar Curvature. Assume furthermore that  the  tangent bundle along the fibres $TN\vert B$ has trivial  Pontrjagin classes in $H^*(N, \R)$.  Then the even differential form on the base manifold $B$ given by
$$
 \int_{N\vert B} d^{-1} \left((\what{A}+a_k \LL) (p(TN\vert B, g))  \right) - \what\eta (\Dir_N) - 2a_k \what\eta (D_N^{sgn})
$$
is closed. Its cohomology class $\sigma (N\vert B, g)\in H^{[0]}(B, \R)$ will be called the higher Kreck-Stolz invariant of the fibre bundle.
\end{definition}

Notice that we have  
$$
d \int_{N\vert B} d^{-1} \left((\what{A}+a_k \LL) (p(TN\vert B, g))  \right) = \int_{N\vert B} (\what{A}+a_k \LL) (p(TN\vert B, g)).
$$
Hence that the differential form of the above definition is closed is a consequence of the relations
$$
d\what{\eta} (\Dir_N) = \int_{N\vert B} \what{A} (TN\vert B, g)\text{ and }2\, d\what{\eta} (D_N^{sgn}) = \int_{N\vert B} \LL (TN\vert B, g).
$$

One can prove many properties similar to the ones listed in \cite{KreckStolz} for a single manifold, including the additivity property for bundles having global sections. All these properties with other applications will be expanded elsewhere, and we  only state the following straightforward proposition:

 \begin{proposition}\
 \begin{enumerate}
\item If $F: (N\to B, g) \rightarrow (N'\to B', g')$ is a bundle map which  preserves the spin structure and is an isometry from $(N, g)$ to $(N', g')$, inducing  $f:B\to B'$,   then $\sigma (N\vert B, g)=f^* \sigma (N'\vert B', g')$.
\item $ \sigma (N\vert B, g)$ only depends on the concordance class  of $g$ in the space of metrics with Positive fiberwise Scalar Curvature.
\end{enumerate}
\end{proposition}

\begin{proof}\
The first item is clear. 
%
Assume that $g_0$ and $g_1$ belong to the same concordance class. So,  the fibre bundle $M=N\times [0, 1]\to B$ can be endowed with a metric $g$ with Positive Fiberwise Scalar Curvature, which coincides with the product metric $g_0+dt^2$ in a collar neighborhood of $N\times \{0\}$ and with $g_1+dt^2$ in a collar neighborhood of $N\times \{1\}$. The Bismut-Cheeger formula then applies for the Dirac operator $\Dir_M$ along the fibers with the fiberwise APS boundary condition to imply that the following differential form is exact (see \cite{BismutCheeger2}  with the normalization from there): 
$$
 \what\eta (\Dir_N, g_1) - \what\eta (\Dir_N, g_0) - \int_{M\times [0, 1]} \what{A}  (TM\vert B, g).
$$
The  family APS index theorem for the fiberwise signature operator also implies that the differential form
$$
2 \what\eta (D^{sgn}_N, g_1) - 2 \what\eta (D^{sgn}_N, g_0) - \int_{M\times [0, 1]} \LL  (TM\vert B, g)
$$
is exact. But since the Pontrjagin classes de $TN\vert B$ are trivial and since the polynomial $\what{A} + a_k \LL$ is a sum of decomposible forms, Lemma \ref{Stockes} applies and we deduce that in our simpler case with $M=N\times [0, 1]$ and since the integral along the fibers of closed forms that vanish at the boundary is a topological cohomology class, the differential form
$$
\int_{M\times [0, 1]} (\what{A} + a_k \LL)  (TM\vert B, g) -\left[ \int_{N\vert B} d^{-1} ((\what{A} + a_k \LL)  (TN\vert B, g_1)) - \int_{N\vert B} d^{-1} ((\what{A} + a_k \LL)  (TN\vert B, g_0))\right]
$$
is also exact. Therefore, combining these statements finishes the proof.
\end{proof}

%
%
%
%
%
%
%

%
%
%

\section{The relative GL-index for families}\label{RelativeFamily}

We give in this last section some remarks and a few results about the expected Gromov-Lawson relative index theorem for families, and we also explain some potential developpements for foliations. The main result is the construction of the higher relative index, although the expected formula for its Chern character is clear, its rigorous proof is too long and tedious to be inserted in this survey paper. It will be a consequence of the more general  computation for foliations in \cite{BenameurHeitsch2021}  which uses the so-called functional calculus of superconnections in terms of the Fourier inversion formula. We restrict ourselves to fibre bundles again where the constructions are already highly non trivial and we shall only state the formula. Our main goal here is hence to define the higher relative index in $K$-theory following a construction due to J. Roe \cite{RoeRelative} since it seemed to us  flexible enough to generalizations. Notice that one can as well define directly the Chern character of the index class without defining the $K$-class as for instance in \cite{BismutCheeger1}.  So we have a smooth fibre bundle $\pi:M\to B$ with compact base. 
Any fiberwise generalized  Dirac operator defines a class in the bivariant $K$-group $KK(M, B)$ when the fiberwise metric is uniformly complete. More precisely,  $D$ extends to a regular self-adjoint operator on the corresponding Hilbert module $\maE$ defined as in the previous sections, and any continuous bounded odd function $f:\R\to \C$ gives rise to  an adjointable odd operator  $f(D)$ acting on  $\maE$. In particular, the operator $F=D(I+D^2)^{-1/2}$ is such a regular self-adjoint operator which is odd for the grading and represents the $KK$-cycle. The following is  classical (see for instance \cite{ConnesSkandalis}):

\begin{proposition}
The triple $(\maE, \maM, F)$, where $\maM$ is the multiplication representation of $C_0(M)$, is an even  Kasparov cycle over the pair $(C_0(M), C_0(B))$ which defines a bivariant class $[D]$ in the Kasparov group $KK (M, B)$. 
\end{proposition}

The pairing of this bivariant class with a compactly supported $K$-theory class $[E_1]-[E_2]$ from $K(M)$ yields again a class in $K(B)$ which is the higher version of the relative index class of Gromov-Lawson $\Ind (D_{E_1}-D_{E_2})$. So, the higher relative index class can be defined easily in this case of a single fibration and one then needs to proceed to compute it. In the case $B=\{\star\}$,  one recovers the Gromov-Lawson relative index defined in terms of parametrices, and hence also as  defined by the cut-and-paste method, see  \cite{GromovLawson3}. We can state the following:

\begin{theorem}\ Assume that $D$ is the spin-Dirac operator along the fibers, supposed to be spin and that  $B$ is a compact manifold, then
$$
\Ch \Ind (D_{E_1}-D_{E_2})=\int_{M\vert B} \what{A} (TM\vert B) \Ch ([E_1]-[E_2]).
$$
\end{theorem}

We point out another possible  method to prove this formula, which is more or less equivalent. It  is  temptating to apply the techniques of the local index theorem in non commutative geometry as proved in \cite{ConnesMoscovici} for a single operator. One then simply applies the family JLO formula in terms of the Bismut superconnection as  proved in  \cite{BenameurCarey}, and then mimic  the computation  of \cite{Ponge} for a single operator. 

We now  explain how to extend the notion of higher relative index to  the general case, say with  two fiberwise Dirac operators on two bundles which agree near infinity.

%
%
%

Let us  fix two smooth fibre bundles $\pi_i:M_i\to B_i$, $i=0,1$ and assume as in the previous section again  that the bases $B_0$ and $B_1$ are compact spaces and that the fibers are even dimensional non-compact complete manifolds. We shall indicate when the smoothness property of the bases is needed. We fix smooth fiberwise uniformly complete metrics $g_0$ and $g_1$ as before. Now the relative index data can be introduced as follows. 

Let $D_i$ be  fiberwise generalized Dirac operators acting on the sections of the $\Z_2$-graded Clifford bundles $E_0=E_0^\pm$  and $E_1=E_1^\pm$ respectively over $M_0$ and $M_1$. Suppose that $\varphi: M_0\smallsetminus K_0\rightarrow M_1\smallsetminus K_1$ is a  homeomorphism (commuting with the projections $\pi_i$) with everywhere invertible vertical tangent map, so that $\varphi$ induces a smooth diffeomorphism between the fibers of $M_0\smallsetminus K_0\to B_0$ and $M_1\smallsetminus K_1\to B_1$  that we further assume to be an isometry between any two fibers. We assume as in the Gromov-Lawson setting of the relative index theorem that $K_0\subset M_0$ and $K_1\subset M_1$ are compact subspaces. Since the fibers are non-compact, $\varphi$ induces a homeomorphism   between $B_0$ and $B_1$ and we can assume  w.l.o.g. that $B_0=B_1=B$.  

As in the previous section, we denote by $\maE_i$ the Hilbert $C(B)$-module associated with the hermitian Clifford bundle $E_i$. 
Denote by $\maE_{K_i}$ and $\maE_{M_i\smallsetminus K_i}$ the submodules of $\maE_i$ composed of all elements of $\maE_i$ which  vanish on $M_i\smallsetminus K_i$ and $K_i$, respectively. Recall that $B$ is compact, we shall assume that  the compact subspaces $K_0$ and $K_1$ are chosen so that the Hilbert submodules  $\maE_{K_0}$ and $\maE_{K_1}$ are unitarily equivalent, i.e. that the associated continuous fields of Hilbert spaces, in the sense of Dixmier, are unitarily isomorphic. 

\begin{remark}
If one is willing to apply the cut-and-paste approach of Gromov and Lawson \cite{GromovLawson3} to define the relative index in this case, then some geometric more precise assumption on $K_0$ and $K_1$ is needed. Namely, one assumes for instance that $K_0$ and $K_1$ are the closures of open subspaces $U_0$ and $U_1$ of $M_0$ and $M_1$ such that  the restriction of $\pi_i$ to $U_i$ is a locally trivial submersion. This condition is satisfied for the family APS problem as well as for many similar index problems which occur in the applications. To fulfil such condition when the bases are compact manifolds, one needs to use a big ball in the ambiant manifold $M_0$ and the corresponding open submanifold in $M_1$. We point out that in this case, $\maE_{K_0}$ and $\maE_{K_1}$ are indeed unitarily equivalent.
\end{remark}

Notice that the Hilbert submodules $\maE_{K_i}$ and $\maE_{M_i\smallsetminus K_i}$  are  full and countably generated. That $\maE_{K_i}$ is for instance  full  is due to our assumption that each fiber of $K_i\to B$ has  positive fiberwise measure. Notice also that we have an orthogonal decomposition $\maE_i=\maE_{K_i}\oplus \maE_{M_i\smallsetminus K_i}$, obtained using multiplication by  the characteristic function of $K_i$ which yields a self-adjoint projection. 

Next we suppose  that we are given a bundle  unitary $\Phi: E_0\vert_{M_0\smallsetminus K_0} \rightarrow E_1\vert_{M_1\smallsetminus K_1}$ over $\varphi$, so continuous and fiberwise smooth in an obvious sense, such that the well defined restrictions of the operators $D_0$ and $D_1$ to $M_0\smallsetminus K_0$ and $M_1\smallsetminus K_1$ are $\Phi$-conjugated.  This is the usual assumption for the relative index to exist. 
Let us then choose a unitary $U$ between $\maE_0$ and $\maE_1$ such that $U$ respects the above orthogonal decompositions associated with the compacts subspaces $K_i$ and such that the induced operator  $\maE_{M_0\smallsetminus K_0}\to \maE_{M_1\smallsetminus K_1}$ coincides with the functorial action given by the identification $(\Phi, \varphi)$.

\begin{definition}\
\begin{itemize}
\item  We denote by $D^*(M_i\vert B, \maE_i)$ the   closure in  $\End_{C(B)} (\maE_i)$ of the subalgebra of adjointable operators with uniform finite propagation.
\item The ideal $C^*(M_i\vert B, \maE_i)$ is composed of those  operators $T\in D^*(M_i\vert B, \maE_i)$ such that 
$$
T f\in \maK (\maE_i)\text{  for any }f\in C_0(M_i).
$$ 
\end{itemize}
\end{definition}

The reader should notice that the above definition of $D^*(M_i\vert B, \maE)$ is different from the usual dual algebra generated by \underline{pseudolocal} finite propagation operators.  The sentence ''$T$ has uniform finite propagation'' means here that there exists a constant $R>0$ such that for any $f_1, f_2\in C_0(M)$ with $
d(\Supp (f_1), \Supp (f_2)) >R$ one has $f_1 T f_2 =0.$
Also $d(Z_1, Z_2)$ is by definition $\infty$ if $\pi (Z_1)\cap \pi (Z_2)=\emptyset$, it is equal to $\inf_{b\in B} d_{M_b} (Z_1\cap M_b, Z_2\cap M_b)$ otherwise. The distance $d_{M_b}$ is of course the riemannian distance associated with the fixed complete metric $g_b$ on the fiber $M_b$. See also \cite{BenameurRoyI}.

\begin{lemma}
For any $T_1\in D^*(M_1\vert B, \maE_1)$ (resp. $C^*(M_1\vert B, \maE_1)$) the operator $U^*T_1U$ belongs to $D^*(M_0\vert B, \maE_0)$ (resp. to $C^*(M_0\vert B, \maE_0)$).
\end{lemma}

\begin{proof}
The proof of this lemma when $B=\{\star\}$  is given in \cite{RoeRelative} and it generalizes immediately to our case. Indeed, if $T_1$ has uniform fiberwise finite propagation then so is $U^*T U$ since $U$ is the functoriality unitary from $\maE_{M_0\smallsetminus K_0}$ to $\maE_{M_1\smallsetminus K_1}$. Moreover, if $f$ is a compactly supported continuous and fiberwise smooth real valued function on $M_0$, then there exists a continuous fiberwise smooth compactly supported function $g:M_1\to \R$ such that $U f= g U f$. Therefore, $
U^*T_1U f = U^*T g U f \in \maK (\maE_0)$ whenever $T\in C^*(M_1\vert B, \maE_1)$.
\end{proof}

\begin{definition}\cite{RoeRelative}
\begin{itemize}
\item We denote by $\maD$ the subspace of  $D^*(M_0\vert B, \maE_0)\times D^*(M_1\vert B, \maE_1)$ composed of all $(T_0, T_1)$ such that the operator $T_0-U^*T_1U$ is a compact operator on the Hilbert module $\maE_0$.
\item We also set $\maC=\maD\cap \left(C^*(M_0\vert B, \maE_0)\times C^*(M_1\vert B, \maE_1)\right)$.
\end{itemize}
\end{definition}

 As a consequence of the previous lemma,  we deduce that $\maD$ is a unital $C^*$-algebra and that $\maC$ is an involutive two-sided ideal in $\maD$. Moreover, the ideal $\maK(\maE_0)$    composed of compact operators  embeds in $C^*(M_0\vert B, \maE_0)$ and hence in $D^*(M_0\vert B, \maE_0)$. Indeed, any two vectors $\xi, \eta$ from $\maE_0$ can be approximated by a couple of compactly supported continuous and fiberwise smooth sections. This shows that any rank one operator on $\maE_0$ can be approximated in the operator norm by a couple of compactly supported fiberwise smoothing operators and this latter is clearly an element of $C^*(M_0\vert B, \maE_0)$. Moreover, $\maK(\maE_0)$ is identified with an involutive two-sided ideal in $\maC$ by using the map $T_0\mapsto (T_0, 0)$. This ideal of compact operators is just the kernel of the second projection, since 
$$
(T_0, 0)\in \maC\;  \Longleftrightarrow\; T_0\in  \maK(\maE_0).
$$
Hence we end up, with the split short exact sequence 
$$
0 \to\;  \maK(\maE_0)\;  \hookrightarrow \; \maC \longrightarrow \; C^*(M_1\vert B, \maE_1)\;  \to 0.
$$
The splitting is  given by the homomorphism $T_1\mapsto (U^*T_1U, T_1)$. Therefore, we deduce in particular that 
$$
K_0 (\maC) \simeq K_0 (\maK (\maE_0) \oplus K_0 (C^*(M_1\vert B, \maE_1)) \simeq K(B) \oplus K_0(C^*(M_1\vert B, \maE_1)),
$$
where the last isomorphism is given by the Morita equivalence $\maK (\maE_0)\sim C(B)$. 

\medskip

Let us now proceed to define the Gromov-Lawson index class $\Ind (D_0, D_1)$ of the relative pair $(D_0, D_1)$ of fiberwise Dirac operators, an element of $K(B)$.

Let $\phi:\R\to \R$ be  a smooth odd non-decreasing function  on $\R$ such that  \cite{RoeBook}
\begin{itemize}
\item $\phi (\pm\infty)=\pm 1$;
 \item The Fourier transform  of $\phi$ is a compactly supported tempered distribution;
\item The function $1-\phi^2$ has smooth compactly supported Fourier transform.
\end{itemize}
By  \cite{BenameurRoyI}[Theorem 4.10], the operator $\phi({D_i})$ is an odd operator which belongs to the $C^*$-algebra $D^*(M_i\vert B, \maE_i)$. In fact it has uniform finite propagation  and is even pseudolocal in the sense that it  commutes with compactly supported continuous functions up to compact operators, while the operator $\phi ({D_i})^2 - I$ belongs  to the $C^*$-algebra $C^*(M_i\vert B, \maE_i)$. 

\begin{proposition}\ We have:
\begin{itemize}
\item $\left(\phi({D_0}), \phi({D_1})\right)\in \maD$. 
\item $\left(I-\phi^2({D_0}), I - \phi^2({D_1})\right)\in \maC$. 
\end{itemize}
\end{proposition}

\begin{proof}
By applying again \cite{BenameurRoyI}[Theorem 4.10], the only remaining property that needs to be proved is that $\phi (D_0)-U^*\phi(D_1) U$ is a compact operator on the Hilbert module $\maE$. But recall our assumption on $K_i$, since $U$ is given by the functoriality unitary off $K_0$ and $K_1$, it is easy to check that there exists a relatively compact open subfibre bundle $L_0$ containing $K_0$  such that for any $\sigma\in \maE_0$ which is supported off $L_0$, we have $U^*\phi (D_1)U \sigma = \phi (D_0)\sigma$. The local compacity of $\phi(D_i)$ then allows to conclude. 
\end{proof}

The positive part of $\left(\phi({D_0}), \phi({D_1})\right)$ then descends to a unitary after moding out by $\maC$. Hence the $K$-theory connecting map from $K_1(\maD/\maC)$ to $K_0(\maC)$, associated with the short exact sequence
$$
0\to \maC \hookrightarrow \maD \longrightarrow \maD/\maC\to 0,
$$
allows to assign to the $K_1$-class of this unitary, an index class in $K_0(\maC)$. But recall that $K_0(\maC) \simeq K(B) \oplus K_0(C^*(M_1\vert B, \maE_1))$. Following \cite{RoeRelative}, we end up with the following definition for families.

\medskip

\begin{definition}
The higher relative index $\Ind (D_0, D_1)$ of the compatible near infinity pair of generalized Dirac operators $D_0$ and $D_1$ is the image of the index class of the pair $\left(\phi({D_0}), \phi({D_1})\right)$ under the projection $K_0(\maC) \longrightarrow  K(B)$ corresponding, in the identification $K_0(\maC) \simeq K(B) \oplus K_0(C^*(M_1\vert B, \maE_1))$, to the first factor. So we have defined our relative index class:
$$
\Ind (D_0, D_1)\; \in \; K(B).
$$
\end{definition}

In order to compute the Chern character of $\Ind (D_0, D_1)$, one works in the smooth category, so in particular $B_0$ and $B_1$ are smooth compact manifolds. The independence of the choice of the function $\phi$  is then the starting point for the  computation of this index class, or rather of its Chern character in $H^{[0]} (B, \R)$. One actually needs to use  functional calculus of superconnections based on the Fourier inversion formula and the Duhamel principle, this is expanded in \cite{BenameurHeitsch2021} for foliations. In our case,  there are simplifications and the expected theorem can be deduced. Let us state the formula for the fiberwise  spin-Dirac operators with coefficients in hermitian bundles which are compatible near infinity (Recall that $\varphi$ induces a diffeomorphism $\varphi_0:B_0\to B_1$):

\begin{theorem}\
We have
$$
\Ch \Ind (\Dir_0\otimes \nabla^0, \Dir_1\otimes \nabla^1) = \int_{M_0\vert B_0} \what{A}(TM_0\vert B_0) \Ch (\nabla^0)\; - \;\varphi_0^*\int_{M_1\vert B_1} \what{A}(TM_1\vert B)\Ch (\nabla^1).
$$
\end{theorem}
The RHS needs to be explained. Since $B_0$ and $B_1$ are compact manifolds, there exists   open submanifolds  $V_0$ and $V_1$ of  $M_0$ and $M_1$ containing respectively $K_0$ and $K_1$   such that  each submersion $\pi_i:V_i\to B_i$ is a submersion and such that $\varphi (M_0\smallsetminus V_0) = M_1\smallsetminus V_1$. Then the experession 
$$
\int_{M_0\vert B_0} \what{A}(TM_0\vert B_0) \Ch (\nabla^0)\; - \;\varphi_0^*\int_{M_1\vert B_1} \what{A}(TM_1\vert B)\Ch (\nabla^1)
$$
 is by definition given by $\int_{V_0\vert B_0} \what{A}(TM_0\vert B_0) \Ch (\nabla^0) - \varphi_0^*\int_{V_1\vert B_1} \what{A} (TM_1\vert B)\Ch (\nabla^1)$ since this difference is well defined and does not depend on the choice of such pair $(V_0, V_1)$.

For general complete foliations, when there exists a holonomy invariant transverse measure, the  relative measured index theorem  has been carried out  by a direct method in \cite{BenameurHeitsch21}.   For the higher version alluded to in the beginning of this section, there are some new complications when using the $K$-index and one rather uses an independently defined higher index in the appropriate cohomology.  Indeed, already in the case of a closed foliated manifold, the cohomological formula, for instance using cyclic homology, is not an obvious corollary of the Connes-Skandalis index theorem in $K$-theory \cite{ConnesSkandalis}, see for instance \cite{ConnesFundamental, BenameurHeitschIa, BenameurHeitschIb, BenameurHeitschII}. In our approach to the higher relative index for general pairs of complete foliations in \cite{BenameurHeitsch2021}, we rather used  a cohomology that is adapted to integration over the leaves, say Haefliger cohomology  \cite{Haefliger, HeitschLazarov}, and the Heitsch superconnection \cite{HeitschSuper} which is the precise  generalization of the Bismut superconnection to foliations built up using the Bott connection.

\appendix

\section{Enlargeability and foliations, a succint overview}

We add in this appendix  a brief overview of recent results about PSC on foliations  which involve the notion of enlargeability. Recall that the classical relative index theorem, say for a single complete manifold, was  used  by Gromov and Lawson in the simplest case of bounded invertibility of the operator to deduce the following breakthrough result:

\begin{theorem}\cite{GromovLawson3}\label{GLenlargeable}
An enlargeable spin closed manifold cannot hold a  PSC metric. 
\end{theorem}

The main observation towards the proof of such theorem is that if $E$ is  trivial near infinity, with Chern character  concentrated in top degree, say $\Ch (E)=\dim (E) + \Ch_{\dim (M)} (E)$ and such that $\Ch_{\dim (M)} (E)$ is a non trivial class, then the formula for the pairing with $[E]-[M\times \C^{\dim (E)}]$ shows that $\Ind (D_{E}-D^N)$ is non-zero. This is true for any bundle that is the pull back from the sphere of  the Bott bundle using a map with non-trivial degree and which is constant near infinity. The notion of enlargeable manifold was introduced precisely to insure that such maps to the sphere do exist (up to going to some covering above $M$) with the extra condition that the pull-back bundle has curvature as small as we please.  The non vanishing of $\Ch_{\dim (M)} (E)$ together with the condition of small curvature then prevents the existence of a metric of uniform PSC, by an easy argument using the Lichnerowicz formula. 

The following generalization to foliations was recently obtained:

\begin{theorem}\cite{BenameurHeitschInv}
A closed manifold $M$ does not admit any spin foliation $F$ with Hausdorff homotopy graph and with a metric whose leafwise scalar curvature satisfies the uniform estimate $> 4/ K_{csf}(M)$. 
\end{theorem}
\medskip

Here the spin condition is on the foliation as in Connes' theorem \cite{ConnesFundamental}, and  $K_{csf}(M)$ is the covering spin foliated $K$-area of $M$, a  variant of the Gromov $K$-area as introduced in \cite{Gromov96} taking the supremum over all spin foliations. Since $K_{csf}(M)=\infty$ when $M$ is enlargeable, a corollary is that a spin foliation, with Hausdorff homotopy graph,  of any enlargeable  closed manifold cannot hold  a metric with  uniform Positive Leafwise Scalar Curvature, see \cite{BenameurHeitschInv}. This obstruction for enlargeable manifolds,  was also obtained in \cite{ZhangEnlargeable}
where Zhang could get rid of the Hausdorff condition, and moreover he was also  able to switch   the spin condition to $M$. 

 There is a similar theorem for complete foliated manifolds excluding Uniformly Positive Leafwise Scalar Curvature \cite{BenameurHeitschExp}.  In \cite{SZ}, Su and Zhang also extended the Zhang results to the non-compact case. 
A theorem for complete manifolds excluding also Positive Leafwise Scalar Curvature could be obtained in  \cite{BenameurHeitschExp}, under the spectral condition that the Novikov-Shubin invariants be  $> 3\codim (F)$ or only $>\codim (F)/2$ for riemannian foliations. 

\bibliographystyle{alpha}

\end{document}